\input ./style/arxiv-general.cfg
\documentclass[MSNbibl,number,citesort,seceqn,dvips]{arxbj}
\makeatletter
   \@ifpackageloaded{graphicx}{}{\usepackage{graphicx}}
\makeatother
\usepackage{mathbh}

\volume{22}
\issue{4}
\pubyear{2016}
\firstpage{2143}
\lastpage{2176}
\doi{10.3150/15-BEJ724}
\docsubty{FLA}

\makeatletter
\newcommand{\rrvert}{\vert}
\newcommand{\llvert}{\vert}
\renewcommand{\mid}{|}
\newtheorem{theorem}{Theorem}[section]
\newtheorem{lemma}[theorem]{Lemma}
\newtheorem{prop}[theorem]{Proposition}
\newremark{remark}[theorem]{Remark}
\newproclaim{assumption}[theorem]{Assumption}
\newcommand{\eqref}[1]{(\ref{#1})}
\newcommand{\E}{\mathbb E}
\newcommand{\V}{\mathbb V}
\newcommand{\Q}{\mathbb Q}
\newcommand{\R}{\mathbb R}
\newcommand{\N}{\mathbb N}
\newcommand{\Z}{\mathbb Z}
\renewcommand{\P}{\mathbb P}

\def\1{\mathbh{1}} 
\newcommand{\ssup}[1]{{(#1)}}
\renewcommand{\d}{{\rm d}}
\newcommand{\eps}{\varepsilon}
\newcommand{\dist}{\operatorname{dist}}
\newcommand{\diam}{\operatorname{diam}}
\newcommand{\Bcal} {{\mathcal B }}
\newcommand{\Ccal} {{\mathcal C }}
\newcommand{\Dcal} {{\mathcal D }}
\newcommand{\Fcal} {{\mathcal F }}
\newcommand{\Mcal} {{\mathcal M }}
\newcommand{\Tcal} {{\mathcal T }}
\newcommand{\Vcal} {{\mathcal V }}
\newcommand{\Wcal} {{\mathcal W }}
\makeatother

\begin{document}
\begin{frontmatter}

\title{Connection times in large ad-hoc mobile~networks}
\runtitle{Connection times in large ad-hoc mobile networks}

\begin{aug}
\author[A]{\inits{H.}\fnms{Hanna}~\snm{D\"oring}\thanksref{A}\ead[label=e1]{hanna.doering@uos.de}},
\author[B]{\inits{G.}\fnms{Gabriel}~\snm{Faraud}\corref{}\thanksref{B}\ead[label=e2]{gabriel.faraud@u-paris10.fr}}
\and
\author[C,D]{\inits{W.}\fnms{Wolfgang}~\snm{K{\"o}nig}\thanksref{C,D}\ead[label=e3]{koenig@wias-berlin.de}}
\address[A]{Universit\"at Osnabr\"uck, Institut f\"ur Mathematik, Albrechtstr. 28a, 49076 Osnabr\"uck, Germany.\\ \printead{e1}}
\address[B]{Laboratoire Modal'x, Universit\'e Paris 10 Nanterre-La D\'
efense, 200 Av. de la R\'epublique, 92000 Nanterre, France. \printead{e2}}
\address[C]{Weierstrass Institute Berlin, Mohrenstr.~39, 10117 Berlin, Germany.\\
\printead{e3}}
\address[D]{Technische Universit\"at Berlin, Institut f\"ur
Mathematik, Str.~des 17.~Juni 136, 10623 Berlin, Germany}
\end{aug}

%
\received{\smonth{1} \syear{2015}}
%
\revised{\smonth{3} \syear{2015}}

%
\begin{abstract}
We study connectivity properties in a probabilistic model for a large
mobile ad-hoc network.
We consider a large number of participants of the system moving
randomly, independently and
identically distributed in a large domain, with a space-dependent
population density of
finite, positive order and with a fixed time horizon. Messages are instantly
transmitted according to a relay principle, that is, they are
iteratively forwarded
from participant to participant over distances smaller than the
communication radius
until they reach the recipient. In mathematical terms, this is a
dynamic continuum percolation model.

We consider the connection time of two sample participants, the amount
of time over
which these two are connected with each other. In the above
thermodynamic limit, we
find that the connectivity induced by the system can be described in
terms of the
counterplay of a local, random and a global, deterministic mechanism,
and we give a formula for the limiting behaviour.

A prime example of the movement schemes that we consider is the
well-known random
waypoint model. Here, we give a negative upper bound for the decay
rate, in the limit
of large time horizons, of the probability of the event that the
portion of the connection time is less than the expectation.
\end{abstract}

%
\end{frontmatter}

\section{Introduction and main results}

\subsection{Background and goals}\label{sec-background}

\textit{Ad-hoc networks} consist of individuals in a given
domain that communicate with each other via a relay principle: messages
are forwarded from individual to individual as long as this
transmission is local, until the message finally arrives at the
recipient. This requires of course that the sender is \textit{connected}
with the recipient, that is, there is a chain of individuals connecting
them such that all links are not larger than a given radius, the \textit{transmission radius} or
\textit{communication radius}. This principle of
message transmission within the system of participants, rather than via
antennas or fixed wires, has a number of advantages over a firmly
installed communication system; for example, its installation is cheap,
it does not require much maintaining, it can accommodate more
information, etc. A disadvantage is of course that the connectivity is
not always fulfilled, that is, it may be that two given individuals are
not connected with each other and are therefore not able to exchange
messages.

The advantages of such a type of system increase if the ad-hoc network
becomes \textit{mobile}, that is, if all the individuals independently
move around in the given region and transmit the messages at their
present location, since in this case a fixed system of wires would be
useless, and firmly located antennas would be necessary, and this may
easily lead to situations of overloads in peak times. This is why
mobile ad-hoc networks are increasingly in the discussion for various
applications, like telecommunication, car-to-car applications for the
distribution of information about the traffic situation, downloading of
large data packages and more \cite{CBD02,Clementi2009,Roy11}. However,
before one can seriously think about an introduction of a mobile ad-hoc
system, one needs to know how reliable it is and how much information
it can reliably transmit and how well the participants of the system
are connected.

The mathematical analysis of the connectivity properties of a
\textit{mobile ad-hoc network} (usually referred to in the literature as
\textit{MANET}), is the purpose of the present paper. We discuss a natural
probabilistic model and derive rigorous results about the quality of
the connection in this system. Roughly speaking, in our model, a large
number $N$ of participants randomly and independently move around in a
given domain $D\subset\R^d$ with $d\geq2$. The movement scheme
considered is quite general, but later we will discuss the prime
example, the \textit{random waypoint model} ({RWP}), in detail. Each of the
participants carries a device that possesses a fixed \textit{communication
radius} $2R$ (the same for everybody). The domain is so large that the
individuals are distributed according to a spatial density that is of
finite order, but may depend on the details of the domain (this models
subareas with more or less frequent visits, like forests, lakes or
public places). We assume that messages are transmitted
instantly, that is, without loss of time. Then we ask, for two fixed
given participants, how large, during a given time interval $[0,T]$,
the amount of time is during which they are connected, their \textit{connection time} $\tau^{\ssup N}_T$. This is one of the most decisive
quantities in such a system, since it measures the quality of the
entire system by means of two sample participants.

The regime in which we will be working is the limit of a large number
of participants, coupled with the limit of a large region such that the
population density (number of participants per area unit) is of finite
positive order. In the language of statistical mechanics, this is the
\textit{thermodynamic limit}. We will condition on the two sample
trajectories. The connection time is obviously a complex function of
the entire system, but we will be able to quantify the influence of the
large number of the participants on the connectivity of the two sample
participants in terms of a simple function. This function is known from
the theory of \textit{continuum percolation}, which studies connectivity
through a union of randomly distributed balls. It turns out that the
limiting connection time has a global, deterministic part and a local,
probabilistic part, the latter of which is described in terms of the
mentioned function. Furthermore, it also turns out that this limit is
deterministic, given the two sample
participants. This is due to one of our assumptions on the movement
scheme, which requires that knowledge about the walker's location at a
later time point does not fix the current location with positive
probability. This assumption implies a certain independence of the
locations of the totality of the walkers at any two given times and
leads to a deterministic limit. This is presented in Sections~\ref
{sec-model} and~\ref{sec-discussion}.

From the practical point of view, a very large value of the connection
time is highly desirable. This can be guaranteed by a large value of
the communication radius $2R$. However, one also would like to have
rules at hand that tell how large this radius must be picked in order
that the connection time exceeds a certain threshold. Some general
answer to this question is given in Section~\ref{sec-discussion}. We
explain there that, under natural conditions, the main effect that may
damage the connection are time lags that any of the two sample
participants spend close to the boundary of the domain $D$, while, in
the interior of $D$, the local connection quality of the system
super-exponentially fast tends to the optimum for $R\to\infty$,
depending on the local user density only.

Furthermore, another important question that we address is about the
long-time behaviour of the connection time. More precisely, we identify
the limiting fraction of the connection time by means of an ergodic
theorem and estimate the probability of the unwanted event that the
connection time covers only an untypically low portion of the time
interval. This is an event of a downward \textit{large deviation}, and we
will show that its probability decays exponentially fast as $T\to
\infty$, and we quantify an upper bound of the decay rate. For this
question, we restrict to the RWP and derive some recurrence properties
that may be useful also for further investigations; see Section~\ref
{sec-RWPResults}.

The model that we consider is sometimes called a dynamic geometric
random graph. Such models were analysed in a series of papers by Peres
and co-workers; see \cite{Peres13}, for example. However, in contrast
to our setting, they do not consider the thermodynamic limit, study
different questions related to the large-time limit, and take Brownian
motion or random walks as the underlying movement scheme. Our purpose
is to study a more realistic model for the random movement of people.

\subsection{Connection time of two participants in the thermodynamic
limit}\label{sec-model}

  Let us introduce the model; our main result here is
Theorem~\ref{thm-connectiontime}.

We consider a system of $N$ particles (the participants of the mobile
ad-hoc network), which randomly move with time horizon $[0,T]$ within a
given bounded open domain $D$ in $\R^d$. The $N$ movements are
independent and identically distributed, and the common movement scheme
(path distribution) does not have to be Markovian; more precise
assumptions follow below. The underlying probability measure and
expectation are denoted by $\P$ and $\E$.

Later (see Section~\ref{sec-RWPResults}) we will be mostly interested
in a particular movement scheme, namely the \textit{random waypoint model}
({RWP}). This motion dynamic is considered in information science as a
realistic model for the random movement of a human being, for example,
a participant of a telecommunication system \cite
{Roy11,Boudec04,BettstetterHartensteinCosta04,BoudecVojnovic06}. A
brief definition of the model is as follows. The walker starts from
some point, picks a random velocity and a random site (the first
waypoint) and then moves with this constant velocity on a straight line
to that waypoint. Having arrived there, he picks the next random
velocity and random waypoint and moves there on a straight line with
the second velocity. This is iterated. All the waypoints are
independent and identically distributed, and the velocities as well,
and the waypoints are independent from the velocities. This model is a
natural extension of the classical RWP, as we admit general
distributions of the waypoints
and the velocities. On the other hand, we do not admit pause times that
the walker spends at waypoints.

Let us proceed with a general movement scheme. We equip every walker
with a fixed communication radius $2R\in(0,\infty)$. That is, there
is a direct connection between any two of them if their distance is at
most $2R$. Two of the $N$ participants, located at $x$ and $y$, say,
are (indirectly) connected if and only if there is a sequence
$x_1,\dots,x_m$ of $m$ other participants such that all the distances
between $x_i$ and $x_{i-1}$ are at most $2R$ for any $i=1,\dots,m+1$,
where we put $x_0=x$ and $x_{m+1}=y$. In other words, the $m+2$ balls
around $x_0,\dots,x_{m+1}$ with radius $R$ have pairwise a non-trivial
intersection along the chain $x_0,\dots,x_{m+1}$; in particular, there
is a continuous path from $x$ to $y$ within their union. This is
fulfilled if and only if $x$ and $y$ lie in the same connected
component of the union of the balls of radius $R$ centred at the $N$
participants. In this way, we see that our model is a dynamic continuum
percolation process.

It is our goal to study the thermodynamic limit of this system, that
is, we think of the volume of the domain being of order $N$, the number
of participants, and we assume that the trajectories are coupled with
$N$ in an accordingly rescaled way. That is, the length scale is
$N^{1/d}$, and the density of participants (their number per unit
volume) is of finite positive order. Then it is clear that a rescaled
version of this picture is better suitable for a mathematical analysis.
Hence, we consider instead the equivalent situation of a fixed domain
$D$ and a fixed movement scheme (both not depending on $N$), and we put
the communication radius equal to $2R N^{-1/d}$. We do not rescale the
time interval $[0,T]$ by $N^{1/d}$, as this is a trivial change.

By $X^{\ssup i}=(X_s^{\ssup i})_{s\in[0,T]}$ we denote the (random)
trajectory of the $i$th participant, that is, a random variable taking
values in the set of functions $[0,T]\to D$. We assume that (making the
underlying probability space $\Omega$ explicit) the map $(s,\omega
)\mapsto X_s^{\ssup i}(\omega)$ is measurable from $[0,T]\times\Omega
$ into $D$. Let $B(x,r)$ denote the open ball around $x$ with radius
$r>0$. Then the set
\[
D_s^{\ssup N}= D\cap\bigcup_{i=1}^N
B \bigl(X_s^{\ssup i}, R N^{-1/d} \bigr) %
\]
is the \textit{communication zone} at time $s$. We introduce the notion of
connectivity at time $s$: for $x,y\in D$ we write
%
\begin{equation}
\label{connectN} x \mathop{\longleftrightarrow}^{N}_{s} y
\quad\Longleftrightarrow \quad x\mbox{ and } y\mbox{ lie in the
same component of }D_s^{\ssup N}.
\end{equation}
We will use this notion only for $x=X^{\ssup1}_s$ and $y=X^{\ssup
2}_s$. Hence, the two participants $X^{\ssup1},X^{\ssup2}$ are
connected at time $s$ if there is a polygon line from $X^{\ssup1}_s$
to $X^{\ssup2}_s$ consisting of line segments of lengths at most
$2RN^{-1/d}$ with the vertices being the locations of other
participants at time~$s$. Hence, $X^{\ssup1}_s\mathop
{\longleftrightarrow}\limits^{N}_{s} X^{\ssup2}_s$ if and only if these two
can exchange a message at time $s$. Note\vspace*{-2pt} that the indicator function $\1
\{X_s^{\ssup1}(\omega) \mathop{\longleftrightarrow}\limits^{N}_{s}
X_s^{\ssup2}(\omega)\}$ is jointly measurable in $s$ and $\omega$,
since it is a polynomial function of the indicators $\1\{\llvert  X_s^{\ssup
i}(\omega)-X_s^{\ssup j}(\omega)\rrvert  \le RN^{-1/d}\}$ with $i,j\in\{
1,\dots,N\}$, which are jointly measurable in $s$ and $\omega$.

The main object is the \textit{connection time}
%
\begin{equation}
\label{taudef} \tau_T^{\ssup N}:=\bigl\llvert \bigl\{ s\in[0,T]
\colon X_s^{\ssup1} \mathop{\longleftrightarrow}^{N}_{s}
X_s^{\ssup2} \bigr\} \bigr\rrvert =\int_0^T\d s \1\bigl\{X_s^{\ssup1} \mathop{\longleftrightarrow
}^{N}_{s} X_s^{\ssup2}\bigr\},
\end{equation}
the amount of time during which these two participants are connected up
to $T$. By the above mentioned joint measurability of the integrand, is
well-defined and measurable. We will analyse the connection time in the
limit $N\to\infty$.

Let us state our assumptions on the random movement of the $N$ walkers.
We write $\{f>r\}$ for short for the set $\{x\in D\colon f(x)>r\}$ and
use analogous notation for similarly defined sets.

\begin{assumption}[(The movement scheme)]\label{Ass-move} The
distribution of the random path $X=(X_s)_{s\in[0,T]}$ in $D$ satisfies
the following:
\begin{longlist}[(ii)]
\item[(i)] For any $s\in(0,T]$, the distribution of the location
$X_s$ possesses a continuous Lebesgue density $f_s\colon D\to[0,\infty)$.

\item[(ii)] For any $s,\widetilde s\in(0,T]$ satisfying $s<\widetilde
s$, a regular version of the conditional distribution of $X_s$ given
$X_{\widetilde s}$ exists and is non-atomic almost surely.
\end{longlist}
\end{assumption}

Sufficient for Assumption~\ref{Ass-move} is the existence of a jointly
continuous Lebesgue density of $X_s$ and $X_{\widetilde s}$ for any
$0<s<\widetilde s\leq T$. Condition (ii) is needed for the asymptotic
independence of the clusters at time $s$ from the clusters at time
$\widetilde s$; it allows us to neglect those walkers that define both
clusters and to deal only with disjoint sets of participants that form
the two clusters; see the proof of Lemma~\ref{lem-prob-conv}.
The reason why it is stated for $s<\widetilde s$ is that it makes the
proof of Lemma~\ref{lem-prob-conv} simpler to understand. It is,
however, possible to adapt it with the same assumption for
$s>\widetilde s$. We leave the details of this to the reader.

We also remark that the map $(s,x)\mapsto f_s(x)$ is measurable.
Indeed, by measurability of the map $(s,\omega)\mapsto X_s(\omega)$,
the indicator $\1 \{\llvert  X_s(\omega)-x\rrvert  \leq\epsilon\}$ is $(\omega, s,
x)$-measurable for any $\epsilon>0$. Therefore, by Fubini's theorem,
its expectation is $(s, x)$-measurable, and by continuity of $f_s$, we
have $f_s(x)=\lim_{\epsilon\downarrow0} \P(\llvert  X_s-x\rrvert  \leq\epsilon
)/(C\epsilon^d)$, $C$ being the volume of the unit ball in $\R^d$,
that is, $f_s(x)$ is a limit of $(s, x)$-measurable functions.

Note that we do not require the continuity of the trajectories;
regularity is only required for the distributions at fixed times.
Assumption~\ref{Ass-move} is satisfied for many diffusions in $D$ and
also for many continuous-time random walks in $D$. For practical
reasons, we are mainly interested in the random waypoint model; see below.

We need to introduce some standard objects from (static, homogeneous)
continuum percolation; see \cite{MeesterRoy96} and Section~\ref
{sec-contperc} below for general background. Let $(Z_i)_{i\in\N}$ be
a standard Poisson point process on $\R^d$ with intensity $\lambda\in
(0,\infty)$. We define the \textit{percolation probability} $\overline
\theta(\lambda,R)$ as the probability that there is a path from
$B(0,R)$ to infinity that never leaves the set $U_R=\bigcup_{i\in\N}
B(Z_i,R)$. In other words, $\overline\theta(\lambda,R)$ is the
probability that $U_R$ has an unbounded connected component that
intersects $B(0,R)$.

The function $\overline\theta$ will play a crucial role in the
asymptotic description of our model. As we will see below, the number
$\overline\theta(\lambda,R)$ describes, in our spatially rescaled
picture, the probability that, locally, a given participant belongs to
the infinitely large cluster and has therefore connection over a
macroscopic part of the space.

We introduce two notions of (non-random) connectedness in the domain
$D$ as follows. By ``path'' we mean a continuous polygon line in
$D$ with finitely many edges, whose vertices lie in $D\cap\Q^d$ (with
possible exception of the first and last one). For $\diamond\in\{\geq
,>\}$ and $x,y\in D$, we write
\[
x \mathop{\longleftrightarrow}^{\diamond}_{s} y \qquad
\Longleftrightarrow\qquad\mbox{there exists a path from }x\mbox{ to }y\mbox{
within }\bigl\{f_s\diamond\lambda_{\rm c}(R)\bigr\}.
\]
Note that the map $(x,y,s)\mapsto\1\{x \mathop{\longleftrightarrow
}\limits^{\diamond}_{s} y\}$ is measurable, as $(s,x)\mapsto f_s(x)$ is and
the notion of a path involves only countably many operations.

Furthermore, we introduce two versions of a limiting value of $\tau
_T^{\ssup N}$. For $\diamond\in\{\geq,>\}$, define
%
\begin{equation}
\label{tauTdef} \tau_T^{\ssup\diamond}\bigl(X^{\ssup1},
X^{\ssup2}\bigr)=\int_0^T \d s \1 \bigl
\{X^{\ssup1}_s \mathop{\longleftrightarrow}^{\diamond}_{s}
X^{\ssup
2}_s \bigr\} \overline\theta^{\ssup\diamond}
\bigl(f_s\bigl(X^{\ssup
1}_s\bigr),R \bigr)
\overline\theta^{\ssup\diamond} \bigl(f_s\bigl(X^{\ssup
2}_s
\bigr),R \bigr),
\end{equation}
where $\overline\theta^{\ssup>}(\lambda,R)=\overline\theta
(\lambda-,R)=\lim_{s \uparrow\lambda}\overline\theta(s,R)$ and
$\overline\theta^{\ssup\geq}(\lambda,R)=\overline\theta(\lambda
+,R)=\lim_{s \downarrow\lambda}\overline\theta(s,R)$ are the left-
and right-continuous versions of $\overline\theta$. Recall that these
two functions coincide at least everywhere outside the critical value
$\lambda_{\rm c}(R)$. Note that $\tau_T^{\ssup\diamond}(X^{\ssup
1}, X^{\ssup2})$ is well-defined and measurable by the measurability
of all the $\theta$-functions, the joint measurabilities of $f_s(x)$
in $s$ and $x$ and of $X_s^{\ssup i}(\omega)$ in $s$ and $\omega$ and
of $\1\{x \mathop{\longleftrightarrow}\limits^{\diamond}_{s} y\}$ in $x$,
$y$ and $s$.

Our main result is the following.

\begin{theorem}\label{thm-connectiontime}Fix $T>0$ and $R>0$, and
assume that the distributions of the $N$ i.i.d.~random movements
$X^{\ssup1},\dots,X^{\ssup N}$ satisfy Assumption~\ref{Ass-move}.
Then, for $\P$-almost all realisations of $(X^{\ssup1}, X^{\ssup
2})$, in probability with respect to $\P(\cdot\mid X^{\ssup1},
X^{\ssup2})$,
%
\begin{equation}
\label{taulimit} \tau_T^{\ssup>}\bigl(X^{\ssup1},
X^{\ssup2}\bigr)\leq\liminf_{N\to\infty} \tau_T^{\ssup N}
\leq\limsup_{N\to\infty}\tau_T^{\ssup N}\leq\tau
_T^{\ssup\geq}\bigl(X^{\ssup1}, X^{\ssup2}\bigr).
\end{equation}
\end{theorem}

For a proof of Theorem~\ref{thm-connectiontime} see Section~\ref
{sec-Proof1}; for a discussion about whether or not the limit in~\eqref
{taulimit} exists and how it behaves for large $R$; see Section~\ref
{sec-discussion}.

The assertion in \eqref{taulimit} shows that the connectivity of the
medium that is built out of $X^{\ssup1},X^{\ssup2},\dots, X^{\ssup
N}$ is fully determined by just two effects: a global, deterministic
one (expressed by the indicator on the event $\{X^{\ssup1}_s \mathop
{\longleftrightarrow}\limits^{\diamond}_{s} X^{\ssup2}_s \}$ in \eqref
{tauTdef}) and a local, stochastic one (expressed by the two $\theta
$-terms). Indeed, the two walkers at time $s$ are connected if and only if:
\begin{enumerate}
\item[$\bullet$] their positions $x=X_s^{\ssup1}$ and $y=X_s^{\ssup
2}$ are connected by a deterministic path within the supercritical
region, that is, the set $\{f_s\diamond\lambda_{\rm c}(R)\}$ (with
$\diamond=\,\geq$ for an upper bound and $\diamond= \,>$ for a lower
bound) and

\item[$\bullet$] both $x$ and $y$ belong locally to the giant
component of the static continuum percolation process with density
$f_{s}(x)$ and $f_{s}(y)$, respectively, and ball radius $RN^{-1/d}$
(note that these two events are asymptotically independent).
\end{enumerate}

\subsection{Discussion}\label{sec-discussion}

\subsubsection{Does the limit in \texorpdfstring{\protect\eqref{taulimit}}{()} exist?}

  Certainly, one expects that, in many cases, $\tau_T^{\ssup
\geq}$ and $\tau_T^{\ssup>}$ should coincide almost surely and
in~\eqref{taulimit} one should have a limit. This is certainly true under
many additional abstract conditions. However, it is difficult to give a
satisfactory sufficient condition that is both reasonably general and
reasonably explicit and, therefore, we abstained from that. Let us
indicate where the difficulties lie.

In order to ensure coincidence of $\tau_T^{\ssup\geq}$ and $\tau
_T^{\ssup>}$, one needs a condition that ensures that connection
within $\{f_s\geq\lambda_{\rm c}(R)\}$ implies connection within $\{
f_s > \lambda_{\rm c}(R)\}$ (at least for the sites $X_s^{\ssup1}$
and $X_s^{\ssup2}$ for almost all $s$) and another condition that
ensures that the $\theta$-terms in \eqref{tauTdef} coincide for
$\diamond=\, >$ and $\diamond=\, \geq$, at least for almost all $s$.

Some sufficient conditions of the first type are certainly easy to
check in many explicit situations, where the structure of the
connectivity landscape given by the density $f_s$ is easy to control.
In general, difficulties can arise if, for $s$ in some set with
positive Lebesgue measure, some components of $\{f_s> \lambda_{\rm
c}(R)\}$ are separated from each other by a component of $\{f_s=
\lambda_{\rm c}(R)\}$ that has a complicated local structure. In
dimension $d=2$, for example, a line with some fractal structure would
pose such a question. In this case, it is unclear what local properties
of the separation set would imply what connectivity probabilities of
the corresponding percolation process. Finding clear criteria seems to
be an open problem in the study of continuum percolation. We believe
that, for related reasons, one can construct situations in which $\tau
_T^{\ssup\geq}$ and $\tau_T^{\ssup>}$ do not coincide, the limit in
\eqref{taulimit} does not exist or is random.

Sufficient conditions of the second type are, in a sense, much easier
to formulate, as the function $\overline\theta(\cdot, R)$ is known
to be continuous outside the critical point $\lambda_{\rm c}(R)$,
\cite{MeesterRoy96}, Theorem~3.9, and, therefore, only times $s$ have
to be considered such that both $X_s^{\ssup1}$ and $X_s^{\ssup2}$ lie
in the set $\{f_s= \lambda_{\rm c}(R)\}$. In fact, in dimension $d=2$,
continuity is also known in the critical point \cite{MeesterRoy96}, Theorem~4.5, such that here the $\theta$-terms do coincide for
any $s$. But in general dimension, continuity in the critical point is
unknown. Hence, in cases where the set $\{f_s= \lambda_{\rm c}(R)\}$
has a positive Lebesgue measure (which can happen only for countably
many values of $R$), there is a positive probability that one of the
two walkers belongs to its interior for a positive portion of the time,
and then the $\theta$-terms may substantially differ.

\subsubsection{Behaviour of the limit in \texorpdfstring{\protect\eqref{taulimit}}{()} for \texorpdfstring{$R\to\infty$}{Rtoinfty}}
From a practical point of view, installing a MANET makes sense only if
the degree of connectivity in the system can be guaranteed to be
extremely high, at least with high probability. Hence, it is a major
goal to find sufficient conditions for a large value (i.e., close to
$T$) of the communication time. Making the communication radius $R$
large is certainly such a criterion, but it is also important to know
how strongly this parameter influences the connectivity. Based on
Theorem~\ref{thm-connectiontime}, we want to illustrate some partial
answer to this question, that is, we want to comment on the behaviour
of the asymptotic lower bound for the connection time, $\tau_T^{\ssup>}$.

This lower bound consists, for any time $s\in[0,T]$, of two
components: the values of $\overline\theta$ in the two locations of
the sample trajectories, and the decision whether or not they are
globally connected through the super-critical area $\{f_s>R^{-d}\lambda
_{\rm c}(1)\}$. An important fact (see \cite{Penrose91}, Corollary of
Theorem 3) is that $\overline\theta(\lambda,R)$ converges
super-exponentially quickly toward 1 for $R\to\infty$, more
precisely, for any $\eps>0$ and some $C_\eps>0$,
%
\begin{equation}
\label{eq-expbounds} \overline\theta(\lambda,R)\geq1-{\operatorname e }^{-\lambda R^d
\llvert  B(0,2)\rrvert  (1-\eps)},
\qquad\lambda R^d\geq C_\eps.
\end{equation}
This shows that the ``bad'' event of being not connected at a
given time $s$ does predominantly not come from the $\overline\theta
$-term, but from the non-connectivity, that is, from the indicator on
the counter-event of $\{X_s^{\ssup1}\mathop{\longleftrightarrow
}\limits^{>}_{s}X_s^{\ssup2}\}$. It is a natural assumption that the density
$f_s$ is, for every $s\in[0,T]$, bounded away from zero in most of the
domain $D$, except\vspace*{1pt} possibly close to the boundary of $D$ and that $f_s$
decays polynomially toward the boundary of $D$. Then the difference
$T-\tau_T^{\ssup>}$ can be upper bounded by some polynomially
decaying term, which depends on the time that at least one of the two
walkers spends polynomially close to the boundary, and some term of the
form ${\operatorname e }^{-C R^d}$ for the remaining time. But the time
that one of the walkers spends close to the boundary of $D$ is
polynomially small in $R$ in probability, since the density is small
there. The conclusion is that bad connectivity properties of the system
predominantly come from the time
that the users spend close to the boundary of $D$, at least if the
domain is homogeneously filled with users.

\subsection{Further investigations for the random waypoint model}\label{sec-RWPResults}

  Let us now concentrate on the random waypoint model, which
was introduced at the beginning of Section~\ref{sec-model}. Below we
show that, under suitable conditions, the RWP is amenable to
Theorem~\ref{thm-connectiontime}, and we study the large-$T$ average
of the connection time and long-time deviations from the mean in terms
of large-deviation estimates.

We have to introduce some notation. We assume that the domain $D$ is
compact and convex. Let $(W_i)_{i\in\N}$ be a sequence of i.i.d. points in $D$, drawn from a distribution $\Wcal$ on $D$, the
\textit{waypoint measure}. Furthermore, let $(V_i)_{i\in\N}$ be an i.i.d. sequence of velocities drawn from some distribution $\Vcal$ on
$(0,\infty)$, the \textit{velocity measure}. The walker starts from an
initial location $X_0\in D$, heading with constant initial velocity
$V_1$ toward the waypoint $W_1$ on a straight line. Having arrived at
$W_1$, the walker immediately moves along the straight line from $W_1$
to $W_2$ with velocity $V_2$ and so on.

This is an extension of the classical RWP, as we admit $D$ as any
convex compact domain, $\Wcal$ as any distribution on $D$, and $\Vcal
$ as any distribution on $(0, \infty)$. On the other hand, we do not
admit pause times that the walker spends at waypoints, as this would
destroy the validity of Assumption~\ref{Ass-move}(ii); in fact, also
the statement of Theorem~\ref{thm-connectiontime} would have to be altered.

We denote by $U_n=\llvert  W_{n+1}-W_n\rrvert  /V_{n+1}$ the time that it takes the
walker to go from the $n$th to the $(n+1)$th waypoint. Then
$T_n=U_0+U_1+\cdots+U_{n-1}$ is the time at which the walker arrives
at the $n$th waypoint, $W_{n}$. We put $T_0=0$. Introduce the
time-change $N(t)=\inf\{n\in\N\colon T_n>t\}$, then $W_{N(t)}$ is
the waypoint that the walker is heading to at time $t$, $V_{N(t)}$ is
his current velocity, and $T_{N(t)}-t$ is the time difference after
which he arrives there. The position of the walker at time $t$ is
denoted by $X_t$. Then
%
\begin{equation}
\label{Xtformula} X_t=W_{N(t)}+\frac{W_{N(t)-1}-W_{N(t)}}{\llvert  W_{N(t)-1}-W_{N(t)}\rrvert  }V_{N(t)}
(T_{N(t)}-t).
\end{equation}
We define all these processes as right-continuous. Note that the
location process $X=(X_t)_{t\in[0,\infty)}$ is not Markov, but the process
%
\begin{equation}
\label{Ydef} Y=(Y_t)_{t\in[0,\infty)}= (X_t,W_{N(t)},V_{N(t)}
)_{t\in
[0,\infty)}
\end{equation}
is a continuous-time Markov process on the state space $\Dcal=D\times
D\times[v_-,v_+]$.

We need to assume some regularity. Throughout the paper, we assume that
the waypoint measure $\Wcal$ and the velocity measure $\Vcal$ possess
continuous Lebesgue densities on $D$ and on some interval
$[v_-,v_+]\subset(0,\infty)$, respectively. In particular, the
velocities are bounded away from 0 and from $\infty$.

We now check that we can apply Theorem~\ref{thm-connectiontime} to the RWP.

\begin{lemma}[(The RWP satisfies Assumption~\ref{Ass-move})]\label{lemmaExampleRWP}
We initialise the RWP by drawing $W_0\in D$ and a velocity $V_0$ from
some distributions on $D$, respectively, on $[v_-,v_+]$ having continuous
densities, such that all the random variables $W_0, W_1, V_0$ are
independent, and put $X_0=W_0$ and $X_t$ as in \eqref{Xtformula}. Then
the RWP satisfies Assumption~\ref{Ass-move}.
\end{lemma}

\begin{pf}
We first show that Assumption~\ref{Ass-move}(i) is satisfied. Indeed,
fix $s\in(0,\infty)$ and note that, on the event $\{s\leq T_1\}$,
\[
X_s=X_0+s V_1\frac{W_1-W_0}{\llvert  W_1-W_0\rrvert  }, %
\]
%
which has obviously a continuous density, since $W_0$, $V_1$ and $W_1$
have and are independent. On the event $\{T_j<s\leq T_{j+1}\}$ with
$j\in\N$, we represent
\[
X_s 
=W_j+(s-T_j)V_{j+1}
\frac{W_{j+1}-W_j}{\llvert  W_{j+1}-W_j\rrvert  }, %
\]
which also has a continuous density, since $W_j$, $V_{j+1}$ and
$W_{j+1}$ have and are independent (and $T_j$ is a continuous function
of them).
Hence, $X_s\1\{T_j<s\leq T_{j+1}\}$ has a continuous density. Summing
on $j\in\N_0$, we also see by use of Dini's theorem that also $X_s$
has a continuous density.

Let us now verify Assumption~\ref{Ass-move}(ii).
For any $x\in D$, $\P(X_s =x \mid X_{\widetilde s}=y)=0$ is clear on the
event $\bigcup_{j\in\N}\{s\leq T_j<\widetilde s\}$, since there was
a change of direction between time $s$ and $\widetilde s$. On the
counter-event, $\bigcup_{j\in\N_0}\{T_{j}<s<\widetilde s \leq
T_{j+1}\}$, we have
\begin{eqnarray*}
\P(X_s =x \mid X_{\widetilde s}=y) &=&\P \biggl(V_{j+1}=
\frac{\llvert  X_{\widetilde s}-x\rrvert  }{\widetilde s-s}, \frac{W_{j+1}-W_j}{\llvert  W_{j+1}-W_j\rrvert  }=\frac{X_{\widetilde
s}-x}{\llvert  X_{\widetilde s}-x\rrvert  } \Bigm| X_{\widetilde s}=y
\biggr)
\\
&\leq&\P \biggl(V_{j+1}=\frac{\llvert  y-x\rrvert  }{\widetilde s-s} \Bigm| X_{\widetilde s}=y
\biggr) =0
\end{eqnarray*}
because the speed is independent from the location and has a continuous density.
\end{pf}

\subsubsection{Long-time limit} Let us consider the long-time
behaviour of $\tau_T^{\ssup\diamond}=\tau_T^{\ssup\diamond
}(X^{\ssup1},X^{\ssup2})$ defined in \eqref{tauTdef} for $\diamond
\in\{>,\geq\}$ for the RWP. We will show in Section~\ref
{invariant-distrib} that the RWP is Harris ergodic and in particular
possesses an invariant distribution, toward which it converges as the
time grows to infinity. In particular, the distribution of the location
of the RWP, $X_t$, converges in total variation sense toward a
probability measure $\mu_*$ on $D$, and it has a continuous Lebesgue
density $f_*\colon D\to[0,\infty)$. However, it is not so easy to
deduce convergence of $\frac{1}T\tau_T^{\ssup\diamond}$ from this,
and we are not able to do so in all cases. For $\diamond\in\{>,\geq\}
$, introduce
%
\begin{equation}
\label{p*def} p_*^{\ssup\diamond}=\int_D \mu_*(\d x)\int
_D \mu_*(\d y) \1\bigl\{ x\mathop{\longleftrightarrow}^{\diamond}_{*}y
\bigr\}\overline\theta ^{\ssup\diamond}\bigl(f_*(x),R\bigr)\overline
\theta^{\ssup\diamond
}\bigl(f_*(y),R\bigr)\in[0,1],
\end{equation}
where $\mathop{\longleftrightarrow}\limits^{\diamond}_{*}$ denotes
connectedness within the set $\{f_*\diamond\lambda_{\rm c}(R)\}$.
Then $p_*^{\ssup\diamond}$ is a measure for connectedness of two
independent sites in $D$ drawn from the limiting distribution of $X_t$.
Furthermore, introduce
%
\begin{equation}
\label{taustardef} \tau_T^{\ssup{\diamond,*}}=\int_0^T
\d s \1\bigl\{X^{\ssup1}_s \mathop{\longleftrightarrow}^{\diamond}_{*}
X^{\ssup2}_s \bigr\} \overline\theta^{\ssup{\diamond}} \bigl(f_*
\bigl(X^{\ssup1}_s\bigr),R \bigr) \overline
\theta^{\ssup{\diamond}} \bigl(f_*\bigl(X^{\ssup2}_s\bigr),R
\bigr),
\end{equation}
the special case of $\tau_T^{\ssup{\diamond}}$ for all the random
waypoint walkers starting in the invariant distribution.

\begin{lemma}[(Ergodic limit)]\label{lem-ergodic}
Let $X^{\ssup1}$ and $X^{\ssup2}$ be two independent copies of $X$. Then
for $\diamond\in\{>,\geq\}$,
%
\begin{equation}
\label{ergodiclimit} \lim_{T\to\infty}\frac{1}T
\tau_T^{\ssup{\diamond,*}}\bigl(X^{\ssup1},X^{\ssup2}
\bigr)=p_*^{\ssup\diamond}\qquad\mbox{almost surely and in }L^1(\P),
\end{equation}
\end{lemma}

We will give a proof of this lemma in Section~\ref
{sec-ergodicproofII}; it is based on a time-discrete Markov chain that
is introduced in Section~\ref{sec-LDP}.

\begin{remark} The previous result is stated with the trajectories of
the walkers started from the invariant state. In general, it is not
clear if $\frac{1}T\tau_T^{\ssup\diamond}$ converges toward
$p_*^{\ssup\diamond}$. Indeed, the critical point is the convergence
of $\1\{x\mathop{\longleftrightarrow}\limits^{\diamond}_{s}y\}$ toward $\1\{
x\mathop{\longleftrightarrow}\limits^{\diamond}_{*}y\}$ for $x,y\in D$ as
$s\to\infty$, which is not true in many counter-examples, as one can
easily find. However, one can check that, under the additional
assumption that $f_s\to f_*$ as $s\to\infty$ uniformly in $D$, then,
in probability,
%
\begin{equation}
\label{ergodiclimitsupandinf} \limsup_{T\to\infty}\frac{1}T
\tau_T^{\ssup\geq}\bigl(X^{\ssup
1},X^{\ssup2}\bigr)\leq
p_*^{\ssup\geq}\quad\mbox{and}\quad\liminf_{T\to\infty}
\frac{1}T\tau_T^{\ssup>}\bigl(X^{\ssup1},X^{\ssup2}
\bigr)\geq p_*^{\ssup>}.
\end{equation}
\end{remark}
%

We remark here that, in cases where the limit in \eqref{taulimit}
exists, we expect that the limits $T\to\infty$ and $N\to\infty$ can
also be interchanged without changing the value, that is,
\[
p_*^{\ssup>}=\lim_{N\to\infty}\lim_{T\to\infty}
\frac{1}T\tau _T^{\ssup N}. %
\]
Indeed, in the limit $T\to\infty$, the ergodic theorem leads to the
average connection probability for two out of $N$ i.i.d.~sites drawn
from the invariant distribution, and then the identification of the
limit $N\to\infty$ follows from Theorem~\ref{thm-connectiontime},
applied to the RWP starting in the invariant distribution. We decided
to leave the details of the proof to the reader.

\subsubsection{Large-$T$ deviations}\label{sec-LDP}

In our next result, we describe the downward deviations of $\tau
_T^{\ssup{>,*}}(X^{\ssup1}, X^{\ssup2})$, more precisely, the
probability of the event $\{\tau_T^{\ssup{>,*}} \leq T p\}$ for $p\in
(0,p_*^{\ssup>})$, in the limit $T\to\infty$. This is certainly an
interesting question, since one would like to effectively bound the
probability of the unwanted event of being connected over less than the
average portion in the long-time limit. We show that this probability
decays even exponentially fast, and we give an explicit bound for the
decay rate. Because of \eqref{taulimit}, such a bound for $\tau
_T^{\ssup{>,*}}$ (rather than for $\tau_T^{\ssup{\geq,*}}$) gives a
useful upper deviation bound for $\tau_T^{\ssup N}$. We write $\P_*$
for the probability measure of the RWP if both copies $Y^{\ssup1}$ and
$Y^{\ssup2}$ start from the invariant distribution.

\begin{theorem}\label{thm-longtimebound} For any $p\in(0,p_*^{\ssup>})$,
%
\begin{equation}
\label{longtimebound} \limsup_{T\to\infty}\frac{1}T\log\P_*\bigl(
\tau_T^{\ssup{>,*}} \leq T p\bigr)<0.
\end{equation}
\end{theorem}

The proof of Theorem~\ref{thm-longtimebound} is in Section~\ref
{sec-proofThDeviat}. It describes an explicit upper bound for the
left-hand side of \eqref{longtimebound} in terms of a variational
problem. The main novelty lies in the proof, which describes the
probability in question in terms of an interesting Markov chain with
nice properties, such that the theory of large deviations may be
applied in a standard way. This Markov chain is an object of
independent interest, as it may serve also for other long-time
investigations of the model, as well as for computer simulations.


\section{Proof of Theorem~\texorpdfstring{\protect\ref{thm-connectiontime}}{1.2}}\label
{sec-Proof1}

In this section, we prove our first main result, Theorem~\ref
{thm-connectiontime}. As a preparation, we first summarise in
Section~\ref{sec-contperc} all relevant available information about
continuum percolation. In Section~\ref{sec-limexptau}, we find the
limit of the expectation of the connection time, and in Section~\ref
{sec-finishThm1} we finish the proof.

\subsection{Static continuum percolation}\label{sec-contperc}

  Let us collect some facts from (static) continuum
percolation; see \cite{MeesterRoy96} or \cite{Penrose03}. Throughout
the paper, we assume that $d\geq2$. Let $(Z_i)_{i\in\N}$ be a
Poisson point process in $\R^d$ with intensity $\lambda>0$. Fix a
radius $R>0$ and consider the union $U_R$ of the balls $B(Z_i,R)$ over
$i\in\N$. We say that two sites $x,y\in\R^d$ are \textit{connected} if
they belong to the same connected component of $U_R$. Connected
components of $U_R$ are called \textit{clusters}. By $\Ccal(x)$, we
denote the cluster that contains $x\in\R^d$. The \textit{percolation
probability} $\theta(\lambda,R)$ is defined as the probability that
$\Ccal(0)$ is unbounded, which we phrase that $0$ is connected with
$\infty$. By scaling, $\theta(\lambda,R)=\theta(\lambda R^d,1)$.
Furthermore, it is known that the map $\lambda\mapsto\theta(\lambda
,R)$ is increasing and that there is a $\lambda_{\rm c}(R)>0$ such
that $\theta(\lambda,R)=0$ for $\lambda<\lambda_{\rm c}(R)$
and $\theta(\lambda,R)>0$ for $\lambda>\lambda_{\rm c}(R)$. Another
characterisation of the critical threshold is that $\llvert  \Ccal(0)\rrvert  =\infty
$ with positive probability for $\lambda>\lambda_{\rm c}(R)$ and
$\llvert  \Ccal(0)\rrvert  <\infty$ with probability 1 for $\lambda<\lambda_{\rm
c}(R)$. In the supercritical case, there exists, with probability one,
a unique cluster with infinite Lebesgue measure, which we call $\Ccal
_\infty$. In the subcritical case, there is no cluster with infinite
Lebesgue measure, almost surely, and the random variable $\llvert  \Ccal(0)\rrvert  $
has finite exponential moments. The map $\lambda\mapsto\theta
(\lambda,R)$ is continuous in any point, with a possible exception at
the critical point, $\lambda_{\rm c}(R)$ \cite{Sarkar97}, Theorem~1.1. The continuity at the critical point is an open
question, but is widely conjectured to be true. For numerical
estimations, we refer to \cite{QZ07}.

Actually, it is not $\theta$ that we will work with in our model, for
the following reason. Certainly, the points $Z_i$ play the role of the
locations of the participants in our telecommunication system. It will
turn out that a given participant located at $Z_i$ is well connected
with the main part of the system if $B(Z_i,R)$ has a non-trivial
intersection with $\Ccal_\infty$; it is not necessary that $Z_i$
itself belongs to $\Ccal_\infty$. Hence, we will be working with a
slightly different notion of percolation: define $\overline\theta
(\lambda,R)$ as the probability that the ball $B(0,R)$ is connected
with $\infty$, that is, that there is an unbounded connected component
of $U_R$ which intersects $B(0,R)$ (obviously, above criticality, this
component has to be $\Ccal_\infty$). Obviously, $\theta\leq
\overline\theta$, and $\overline\theta$ shares the above mentioned
properties with $\theta$; however, with possibly different numerical
values. In particular, $\overline\theta$ possesses the same scaling
properties, and is an increasing function of $\lambda$, and is
positive above some threshold and zero below. One can also easily check
that the percolation threshold is the same for the two definitions, and
that the proof of the continuity for the usual definition extends to
this definition.

\subsection{Limiting expectation of the connection time}\label{sec-limexptau}

  We fix $T>0$ for the remainder of the section. In the
following, we abbreviate
\[
\P_{1,2}(\cdot)=\P \bigl(\cdot \mid X^{\ssup1},X^{\ssup2}
\bigr)\quad\mbox{and}\quad\E_{1,2}[ \cdot ]=\E \bigl[\cdot \mid
X^{\ssup1},X^{\ssup2} \bigr]. %
\]
Use \eqref{taudef} and Fubini's theorem to see that
\[
\E_{1,2}\bigl[\tau_T^{\ssup{N}}\bigr]=\int_0^T\d s \P_{1,2}
\bigl(X^{\ssup
1}_s\mathop{\longleftrightarrow}^{N}_{s}
X^{\ssup2}_s \bigr). %
\]
%
We are going to approximate the event $\{X^{\ssup1}_s\mathop
{\longleftrightarrow}\limits^{N}_{s} X^{\ssup2}_s\}$ by the event that
$X^{\ssup1}_s$ and $X^{\ssup2}_s$ are separated from each other, but
connected through either $\{f_s>\lambda_{\rm c}(R)\}$ or through $\{
f_s\geq\lambda_{\rm c}(R)\}$ and belong locally to the macroscopic
part of the communication zone. More precisely, for $s\in[0,T]$,
$\delta>0$ and $N\in\N$, we introduce the events
\[
G_{N,s,\delta}^{\ssup i}=\bigl\{X_s^{\ssup i}\mathop
{\longleftrightarrow}^{N}_{s}\partial
\bigl[X_s^{\ssup i}+(-\delta /2,\delta/2)^d \bigr]
\bigr\}, \qquad i\in\{1,2\}, %
\]
that $X_s^{\ssup i}$ and at least some point of the boundary of the
$\delta/2$-box around $X_s^{\ssup i}$ lie in the same connected
component of the union of the $RN^{-1/d}$-balls around $X_s^{\ssup
1},\dots,X_s^{\ssup N}$.
Note that $G_{N,s,\delta}^{\ssup i}$ is not modified by changing the
locations of the walkers outside the $\delta$-box around $X_s^{\ssup i}$.

We will give bounds for the connection time $\tau_T^{\ssup N}$ in
terms of
\begin{eqnarray*}
\tau_T^{\ssup{N,\delta,\diamond}}
\bigl(X^{\ssup1}, X^{\ssup2}\bigr) &=& \int_0^T
\d s \prod_{i=1}^2 \1_{G_{N,s,\delta}^{\ssup i}} \1
\bigl\{\bigl\llvert X_s^{\ssup1}-X_s^{\ssup2}
\bigr\rrvert \geq3\delta\bigr\}\1\bigl\{X_s^{\ssup
1}\mathop{
\longleftrightarrow}^{\diamond}_{s}X_s^{\ssup2}
\bigr\},
\end{eqnarray*}
%
in the limit $N\to\infty$, followed by $\delta\downarrow0$. We will
use $\tau_T^{\ssup{N,\delta,>}}$ as a lower bound and $\tau
_T^{\ssup{N,\delta,\geq}}$ as an upper bound for $\tau_T^{\ssup
N}$. Recall the quantities $\tau_T^{\ssup\diamond}$ defined in
\eqref{tauTdef}, which will serve as limiting objects of $\tau
_T^{\ssup{N,\delta,\diamond}}$.

\begin{prop}[(Limiting expectation of $\tau_T^{\ssup N}$)]\label
{prop-exptau}Let the distributions of the $N$ i.i.d.~walkers satisfy
Assumption~\textup{\ref{Ass-move}(i)}. Then, for $\P$-almost all $X^{\ssup1}$
and $X^{\ssup2}$, provided that $R$ is chosen such that $\int_0^T \d
s  \1\{f_s(X_s^{\ssup i})=\lambda_{\rm c}(R)\} =0$, for $i=1,2$:
\begin{longlist}[(ii)]
\item[(i)]
%
\begin{eqnarray}
\limsup_{\delta\downarrow0}\limsup_{N\to\infty}
\E_{1,2} \bigl(\tau_T^{\ssup N}-\tau_T^{\ssup{N,\delta,\geq}}
\bigr)^+&=& 0,\label{firsttauesti}
\\
\liminf_{\delta\downarrow0}\liminf_{N\to\infty}
\E_{1,2} \bigl(\tau_T^{\ssup N}-\tau_T^{\ssup{N,\delta,>}}
\bigr)^-&=& 0.\label
{secondtauesti}
\end{eqnarray}

\item[(ii)] For any $\diamond\in\{>,\geq\}$,
%
\begin{equation}
\lim_{\delta\downarrow0}\lim_{N\to\infty} \E_{1,2}
\bigl[\tau _T^{\ssup{N,\delta,\diamond}}\bigr]=\tau_T^{\ssup\diamond}
\bigl(X^{\ssup
1}, X^{\ssup2}\bigr).
\end{equation}
\end{longlist}
\end{prop}

The main step in the proof is the following.

\begin{lemma}\label{lem-Penrose} Let the distributions of the $N$
i.i.d.~walkers satisfy Assumption~\textup{\ref{Ass-move}(i)}. Then, for $\P
$-almost all $X^{\ssup1}$ and $X^{\ssup2}$, for almost any $s\in
[0,T]$ and on the event $\{f_s(X_s^{\ssup1})\neq\lambda_{\rm c}(R)\}
\cap\{f_s(X_s^{\ssup2})\neq\lambda_{\rm c}(R)\}\cap\{X_s^{\ssup
1}\neq X_s^{\ssup2}\}$:
\begin{longlist}[(ii)]
\item[(i)]
%
\begin{eqnarray}
&& \limsup_{\delta\downarrow0}\limsup_{N\to\infty}
\P_{1,2} \bigl[ \bigl(X_s^{\ssup1}\mathop{
\longleftrightarrow}^{N}_{s} X_s^{\ssup 2}
\bigr)\setminus \bigl(G_{N,s,\delta}^{\ssup1}\cap G_{N,s,\delta}^{\ssup2}
\cap\bigl\{X_s^{\ssup1}\mathop {\longleftrightarrow}^{\geq}_{s}X_s^{\ssup2}
\bigr\} \bigr) \bigr]
\nonumber\\[-8pt]\label{eq-limsup}  \\[-8pt]\nonumber
&&\quad = 0,
\\
&& \limsup_{\delta\downarrow0}\limsup_{N\to\infty}
\P_{1,2} \bigl[ \bigl(G_{N,s,\delta}^{\ssup1}\cap
G_{N,s,\delta}^{\ssup2}\cap\bigl\{ X_s^{\ssup1}
\mathop{\longleftrightarrow}^{>}_{s}X_s^{\ssup2}
\bigr\} \bigr)\setminus \bigl(X_s^{\ssup1}\mathop{
\longleftrightarrow }^{N}_{s} X_s^{\ssup2}
\bigr) \bigr]
\nonumber\\[-8pt]\label{eq-liminf} \\[-8pt]\nonumber
&&\quad = 0.
\end{eqnarray}

\item[(ii)]
%
\begin{eqnarray}\label{eq-Penrose}
\overline\theta\bigl(f_s
\bigl(X_s^{\ssup1}\bigr)-,R\bigr)\overline\theta
\bigl(f_s\bigl(X_s^{\ssup2}\bigr)-,R\bigr) &\leq&
\liminf_{\delta\downarrow0}\liminf_{N\to\infty}\P
_{1,2} \bigl(G_{N,s,\delta}^{\ssup1}\cap G_{N,s,\delta}^{\ssup2}
\bigr)\nonumber
\\
&\leq&\limsup_{\delta\downarrow0}\limsup_{N\to\infty}\P
_{1,2} \bigl(G_{N,s,\delta}^{\ssup1}\cap G_{N,s,\delta}^{\ssup2}
\bigr)
\\
&\leq& \overline\theta\bigl(f_s\bigl(X_s^{\ssup1}
\bigr)+,R\bigr)\overline\theta \bigl(f_s\bigl(X_s^{\ssup2}
\bigr)+,R\bigr).\nonumber
\end{eqnarray}
\end{longlist}
\end{lemma}

\begin{pf}
Fix $s$ and let us abbreviate $x=X^{\ssup1}_s$ and $y=X^{\ssup2}_s$.
Under $\P_{1,2}$, only the sites $X^{\ssup3}_s,\dots,X_s^{\ssup N}$
are random (in fact, they are i.i.d.~with density $f_s$), but the
notion of connectedness and components induced by the point process
refer to \textit{all} the balls $B(X^{\ssup i}_s,RN^{-1/d})$ with
$i=1,2,\dots,N$.

Let us prove (ii).
First, we consider the case that $f_s(x)<\lambda_{\rm c}(R)$ or
$f_s(y)< \lambda_{\rm c}(R)$, in which case the events $\{X_s^{\ssup
1}\mathop{\longleftrightarrow}\limits^{>}_{s}X_s^{\ssup2}\}$ and $\{
X_s^{\ssup1}\mathop{\longleftrightarrow}\limits^{\geq}_{s}X_s^{\ssup2}\}$
are not fulfilled. Without loss of generality, let us assume that
$f_s(x)<\lambda_{\rm c}(R)$. Choose $\delta>0$ so small that the
$\delta$-box around $x$ does not contain $y$ and that $f_s<\lambda
_{\rm c}(R)$ within that box. We apply \cite{Penrose95}, Proposition 2,
for $\eps=\delta/4$ and obtain that, with $\P_{1,2}$-probability
tending to $1$ as $N\to\infty$, any connected component of $\bigcup_{i=3}^N B(X^{\ssup i}_s,RN^{-1/d})$ in this cube has a diameter
bounded from above by $\eps$. In particular, with $\P
_{1,2}$-probability tending to $1$, $x$ is not connected with the
boundary of the cube $x+(-\delta,\delta)^d$. 
Therefore, \eqref{eq-Penrose} is trivial, as all terms are zero.

To\vspace*{1pt} prove \eqref{eq-Penrose} in the remaining case $f_s(x)\geq\lambda
_{\rm c}(R)$ and $f_s(y)\geq\lambda_{\rm c}(R)$, we show now that the
two events $G_{N,s,\delta}^{\ssup1}$ and $G_{N,s,\delta}^{\ssup2}$
are asymptotically independent with $\P_{1,2}$-probabilities tending
to $\overline\theta(f_s(x),R)$ and $\overline\theta(f_s(y),R)$,
respectively. Let $\mu_s$ denote the measure with density $f_s$.
Indeed, first note that, for every sufficiently large $N$ such that the
ball diameter $2RN^{-1/d}$ is less than the distance between
$x+(-\delta,\delta)^d$ and $y+(-\delta,\delta)^d$. Hence, the
positions of the points falling in $x+(-\delta,\delta)^d$ and
$y+(-\delta,\delta)^d$ are independent, conditionally on their
numbers. These two numbers are binomially distributed with parameters
$N$ and $\mu_s(x+(-\delta,\delta)^d))$ and $\mu_s(y+(-\delta
,\delta)^d)$, respectively. Therefore, by the law of large numbers,
they stochastically dominate, with \mbox{$\P_{1,2}$-}probability tending to
$1$, the Poisson law with parameters $N(\mu_s(x+(-\delta,\delta
)^d)-\eta(2\delta)^d)$ and $N(\mu_s(y+(-\delta,\delta)^d)-\eta
(2\delta)^d)$, respectively, for any $\eta>0$. Note that the events
$G_{N,s,\delta}^{\ssup
1}$ and $G_{N,s,\delta}^{\ssup2}$ are monotonic\vspace*{1pt} in the intensity,
that is, their $\P_{1,2}$-probability is not larger than the $\P
_{1,2}$-probability of the same event under continuum percolation in
$x+(-\delta,\delta)^d$ and $y+(-\delta,\delta)^d$ with intensity
parameters $f_s(x)-2\eta$ and $f_s(y)-2\eta$, respectively, and ball
diameter $RN^{-1/d}$. Since we are now considering Poisson point
processes, the events are independent. Their respective probabilities
converge toward $\overline\theta(f_s(x)-2\eta,R)$ and $\overline
\theta(f_s(y)-2\eta,R)$. Since this is true for any $\eta$, we can
use the continuity of $ \overline\theta(\cdot,R)$, to obtain the
lower bound in \eqref{eq-Penrose}. The upper bound is proved in a
similar manner, using that $\overline\theta(\lambda)$ is the
limiting probability that the origin is connected with the boundary of
a centred cube for diverging radius. This finishes the proof of (ii).

In order to show (i), we are going to decompose into four separate cases.
First, we consider the case that $f_s(x)<\lambda_{\rm c}(R)$ or
$f_s(y)< \lambda_{\rm c}(R)$. As before, let us assume that
$f_s(x)<\lambda_{\rm c}(R)$. With $\P_{1,2}$-probability tending to
$1$, $x$ is not connected with the boundary of the cube $x+(-\delta
,\delta)^d$ and, therefore, neither with $y$, by the previous
argument. This proves \eqref{eq-limsup} and \eqref{eq-liminf} in this case.

In the second part of the proof, we assume that $x$ and $y$ belong to
the same component of $\{f_s> \lambda_{\rm c}(R)\}$, in which case
both events $\{X_s^{\ssup1}\mathop{\longleftrightarrow
}\limits^{>}_{s}X_s^{\ssup2}\}$ and $\{X_s^{\ssup1}\mathop
{\longleftrightarrow}\limits^{\geq}_{s}X_s^{\ssup2}\}$ are fulfilled. Pick
some auxiliary parameter $\eta>0$ that is smaller than $f_s(x)-\lambda
_{\rm c}(R)$ and smaller than $f_s(y)-\lambda_{\rm c}(R)$. Now, using
the continuity of $f_s$ in accordance with Assumption~\ref
{Ass-move}(i), pick $\delta>0$ so small that $x+(-\delta,\delta)^d$
and $y+(-\delta,\delta)^d$ have positive distance and that $f_s$
takes values in $[f_s(x)-\eta,f_s(x)+\eta]$ in $x+(-\delta,\delta
)^d$ and values in $[f_s(y)-\eta,f_s(y)+\eta]$ in $y+(-\delta,\delta
)^d$ and such that there exists a set of the form $U=\bigcup_{i=0}^m
2\delta z_i+[-\delta,\delta]^d$ in $\{f_s> \lambda_{\rm c}(R)\}$
with $m\in\N$, $z_1,\dots, z_m\in\Z^d$ such that $z_i$ and
$z_{i-1}$ are nearest neighbours for any $i=1,\dots, m$ and $x+(-
\delta,\delta)^d\subset U$ and $y+(-\delta,\delta)^d\subset U$ and
$f_s>\lambda_{\rm c}(R)$ inside $U$. That this is possible is
easy to see by elementary continuity and compactness arguments. Since
$U$ is a compact subset of $\{f_s> \lambda_{\rm c}(R)\}$, the density
$f_s$ is even bounded away from $\lambda_{\rm c}(R)$ on $U$.

Let $\Ccal_{x,\delta}^{\ssup{s,N}}$ and $\Ccal_{y,\delta}^{\ssup{s,N}}$, respectively, denote the largest component of the union of
the $RN^{-1/d}$-balls around the points $X_s^{\ssup1},\dots
,X_s^{\ssup N}$ which lie in $x+(-\delta,\delta)^d$, respectively, in
$y+(-\delta,\delta)^d$. According to \cite{Penrose95}, Proposition~3, with $\P_{1,2}$-probability tending to $1$ as $N\to
\infty$, these are the only ones in the respective boxes whose size
(measured in terms of the number of $i$ such that $X_s^{\ssup i}$
belongs to it) is of order $N$, and they are also uniquely determined
by requiring their diameter of positive order. In particular, as $N\to
\infty$, the probability of the symmetric difference between the
events $\{x \in\Ccal_{x,\delta}^{\ssup{s,N}}\}$ and $G_{N,s,\delta
}^{\ssup1}$ (resp., $\{y\in\Ccal_{y,\delta}^{\ssup{s,N}}\}$
and $G_{N,s,\delta}^{\ssup2}$) goes to zero. By \cite{Penrose95}, Proposition~4, such a unique cluster, $\Ccal_U^\ssup{s,N}$ also
exists for the set $U$. Hence,
with $\P_{1,2}$-probability tending to $1$, both $\Ccal_{x,\delta
}^\ssup{s,N}$ and $\Ccal_{y,\delta}^\ssup{s,N}$ belong to $\Ccal
_U^\ssup{s,N}$. This implies that with probability tending to $1$ as
$N\to\infty$, the symmetric difference between the event $\{x\mathop
{\longleftrightarrow}\limits^{N}_{s} y\}$ and the event $G_{N,s,\delta
}^{\ssup1}\cap G_{N,s,\delta}^{\ssup2}$ goes to zero, which implies
\eqref{eq-liminf} and \eqref{eq-limsup}.

In the third case, we have $f_s(x)>\lambda_{\rm c}(R)$ and $f_s(y)>
\lambda_{\rm c}(R)$, and $x\mathop{\longleftrightarrow}\limits^{\geq}_{s}
y$, but not $x\mathop{\longleftrightarrow}\limits^{>}_{s} y$, in which case
\eqref{eq-liminf} is trivial, as the event inside the probability is empty.
To prove \eqref{eq-limsup}, it is enough to see that,
deterministically, the existence of a path between $x$ and $y$ implies
$G_{N,s,\delta}^{\ssup1}$ and~$G_{N,s,\delta}^{\ssup2}$. %

In the fourth case, we have $f_s(x)>\lambda_{\rm c}(R)$ and $f_s(y)>
\lambda_{\rm c}(R)$, but not $x\mathop{\longleftrightarrow}\limits^{\geq
}_{s} y$. Here, \eqref{eq-liminf} is again trivial, as the event
inside the probability is empty. To prove \eqref{eq-limsup}, it is
enough to check that, with probability tending to $1$, $x$ and $y$ are
not connected in the union of the $RN^{-1/d}$-balls around the points
$X_s^{\ssup1},\dots,X_s^{\ssup N}$. Here, it is intuitively clear
that any path between $x$ and $y$ has to cross a non-trivial zone where
$f_s<\lambda_{\rm c}(R)$ and that this disconnects $x$ and $y$ in the
limit. Let us give a proof.

First, we argue that there is a (deterministic) compact set $\Gamma
\subset D$ and $\eps,\gamma>0$ such that $\Gamma\subset\{f_s\leq
\lambda_{\rm c}(R)-\eps\}$ and every path connecting $x$ and $y$
passes through $\Gamma$ for at least $ \gamma$ space units. Indeed,
since $x\mathop{\longleftrightarrow}\limits^{\geq}_{s} y$ does not hold,
$x$ and $y$ lie in disjoint components of $\{f_s\geq\lambda_{\rm
c}(R)\}$. Hence, both these components have a positive distance $\eta$
to the remainder of $\{f_s\geq\lambda_{\rm c}(R)\}$, since these
three sets are compact and mutually disjoint. Abbreviate
\[
\Gamma_\alpha=\bigl\{z\in D\colon\dist\bigl(z,\bigl\{f_s
\geq\lambda_{\rm c}(R)\bigr\} \bigr)\geq\alpha\bigr\},\qquad\alpha>0,
\]
and pick $\Gamma=\Gamma_{\eta/16}$. Then every path from $x$ to $y$
passes at least a distance $\gamma=\eta-2\eta/16=7\eta/8$ through
$\Gamma$. By continuity of $f_s$, this set $\Gamma$ is compact and is
contained in $\{f_s\leq\lambda_{\rm c}(R)-\eps\}$ for some $\eps>0$.

Second, we argue that, with $\P_{1,2}$-probability tending to 1 as
$N\to\infty$, any connected component of $\bigcup_{i=3}^N B(X^{\ssup
i}_s,RN^{-1/d})$ in $\Gamma$ has\vspace*{2pt} a diameter at most $\gamma/2$.
Indeed, consider the neighbourhood $\widetilde\Gamma=\Gamma_{\eta
/32}$ of $\Gamma$, then, for $N$ sufficiently large, the connected
components inside $\Gamma$ do not depend on the configuration outside
$\widetilde\Gamma$. By continuity of $f_s$, on $\widetilde\Gamma$,
the function $f_s$ is still bounded away from $\lambda_{\rm c}(R)$,
say it is bounded from above by $\lambda_{\rm c}(R)-\widetilde\eps$
for some $\widetilde\eps>0$. We upper bound the probability of having
any connected component inside $\widetilde\Gamma$ of diameter bigger
than $\gamma/2$ against the same probability under the homogeneous
Poisson point process with intensity parameter $\lambda_{\rm
c}(R)-\widetilde\eps/2$ on some cube that contains $\widetilde\Gamma
$ (see the above argument). Now, as this intensity parameter is subcritical,
this
probability tends to $0$ as $N\to\infty$.

Now we finish the proof of \eqref{eq-limsup} and \eqref{eq-liminf} in
the fourth case. Indeed, the existence of a connection from $x$ to $y$
through $\bigcup_{i=1}^N B(X^{\ssup i}_s,RN^{-1/d})$ implies the
existence of at least one connected component of this set in $\Gamma$
of diameter at least $\gamma$, since any path from $x$ to $y$ passes
at least a distance $\gamma$ through $\Gamma$. But, as we saw in the
second step, the probability of this existence tends to $0$ as $N\to
\infty$.
\end{pf}

\begin{pf*}{Proof of Proposition~\ref{prop-exptau}}
Observe that
\begin{eqnarray*}
&&\E_{1,2} \bigl(\tau_T^{\ssup N}-
\tau_T^{\ssup{N,\delta,\geq}} \bigr)^+
\\
&&\quad \leq \int_0^T \d s \bigl(\P_{1,2}
\bigl[ \bigl(X_s^{\ssup1}\mathop {\longleftrightarrow}^{N}_{s}
X_s^{\ssup2} \bigr)\setminus \bigl(G_{N,s,\delta}^{\ssup1}
\cap G_{N,s,\delta}^{\ssup2}\cap\bigl\{ X_s^{\ssup1}
\mathop{\longleftrightarrow}^{\geq}_{s}X_s^{\ssup2}
\bigr\} \bigr) \bigr]
\\
&&\qquad{}\times \1 \bigl\{\bigl\llvert X^{\ssup1}_s-X^{\ssup2}_s
\bigr\rrvert > 3\delta\bigr\} \1 \bigl\{f_s\bigl(X_s^{\ssup1}
\bigr)\neq\lambda_{\rm c}(R)\bigr\} \1 \bigl\{f_s
\bigl(X_s^{\ssup2}\bigr)\neq\lambda_{\rm c}(R)\bigr\}
\\
&&\qquad{} + \1 \bigl\{\bigl\llvert X^{\ssup1}_s-X^{\ssup2}_s
\bigr\rrvert < 3\delta\bigr\} + \1 \bigl\{f_s\bigl(X_s^{\ssup1}
\bigr)= \lambda_{\rm c}(R)\bigr\} + \1 \bigl\{f_s
\bigl(X_s^{\ssup2}\bigr)= \lambda_{\rm c}(R)\bigr\}
\bigr).
\end{eqnarray*}
Hence, by \eqref{eq-limsup},
\[
\limsup_{\delta\downarrow0}\limsup_{N\to\infty}
\E_{1,2} \bigl(\tau_T^{\ssup N}-\tau_T^{\ssup{N,\delta,\geq}}
\bigr) \leq\int_0^T \d s \1 \bigl\{\bigl\llvert
X^{\ssup1}_s-X^{\ssup2}_s\bigr\rrvert =0
\bigr\}, %
\]
according to our assumption on $R$. Note that, almost surely, $\int_0^T \d s  \1 \{\llvert  X^{\ssup1}_s-X^{\ssup2}_s\rrvert  =0\} =0$, since $X^{\ssup
1}_s$ and $X^{\ssup2}_s$ are independent with density $f_s$ for any
$s\in[0,T]$. Hence, the proof of \eqref{firsttauesti} is finished.
The proof of \eqref{secondtauesti} is done in the same way using
\eqref{eq-liminf}. Hence, part (i) is proved.

Now we turn to the proof of (ii).

Note that our assumptions exclude that $f_s(X^{\ssup i}_s)=\lambda
_{\rm c}(R)$ outside a set of measure zero. Therefore this does not
appear in the integral. Furthermore, 
$\overline\theta$ is continuous except maybe for $\lambda_{\rm
c}(R)$. Therefore, for almost every $s$, \eqref{eq-Penrose}
reformulates to
\begin{eqnarray*}
\lim_{\delta\downarrow0}\liminf_{N\to\infty}
\P_{1,2} \bigl(G_{N,s,\delta}^{\ssup1}\cap G_{N,s,\delta}^{\ssup2}
\bigr)&=& \lim_{\delta\downarrow0}\limsup_{N\to\infty}
\P_{1,2} \bigl(G_{N,s,\delta}^{\ssup1}\cap G_{N,s,\delta}^{\ssup2}
\bigr)
\\
&=&\overline\theta\bigl(f_s\bigl(X_s^{\ssup1}
\bigr),R\bigr)\overline\theta \bigl(f_s\bigl(X_s^{\ssup2}
\bigr),R\bigr).
\end{eqnarray*}

Thus, (ii) follows by Lebesgue's theorem.
\end{pf*}

\subsection{Finish of the proof}\label{sec-finishThm1}
The second main step in proving Theorem~\ref{thm-connectiontime} is
the following lemma.
Recall that $\P_{1,2}$ denotes the conditional distribution given
$X^{\ssup1}$ and $X^{\ssup2}$.

\begin{lemma}[($\tau_T^{\ssup{N,\delta,\diamond}}$ is asymptotically
deterministic)]\label{lem-prob-conv}Let the distributions of the $N$
i.i.d.~walkers satisfy Assumption~\textup{\ref{Ass-move}(i)} and \textup{(ii)}. Then,
for any $\diamond\in\{>,\geq\}$, for almost every paths $ X^{\ssup
1},X^{\ssup2}$, the difference
$\tau_T^{\ssup{N,\delta,\diamond}}-\E_{1,2}[\tau_T^{\ssup
{N,\delta,\diamond}}]$ vanishes as $N\to\infty$, followed by
$\delta\downarrow0$, in $\P_{1,2}$-probability, provided that $R$ is
chosen such that
$\int_0^T \d s  \1\{f_s(X_s^{\ssup i})=\lambda_{\rm c}(R)\} =0$ for $i=1,2$.
\end{lemma}

\begin{pf}
The claimed convergence follows, by Chebyshev's inequality, from the
fact that the $\P_{1,2}$-variance of $\tau_T^{\ssup{N,\delta,\diamond}}$ vanishes. Writing $\V_{1,2}$ for the $\P
_{1,2}$-variance, this is equal to
%
\begin{eqnarray}\label{variance}
\V_{1,2}\bigl(\tau^{\ssup{N,\delta,\diamond}}_T\bigr) &=& \int_0^T
\d s\int_0^T\d\widetilde s \1\bigl\{\bigl\llvert
X_s^{\ssup1}-X_s^{\ssup2}\bigr\rrvert >3
\delta\bigr\}\1\bigl\{ X_s^{\ssup1}\mathop{
\longleftrightarrow}^{\diamond}_{s}X_s^{\ssup
2}\bigr\}\nonumber
\\
&&{}\times \1\bigl\{\bigl\llvert X_{\widetilde s}^{\ssup1}-X_{\widetilde s}^{\ssup
2}\bigr\rrvert >3\delta\bigr\}\1\bigl\{X_{\widetilde s}^{\ssup1}\mathop {
\longleftrightarrow}^{\diamond}_{\widetilde s}X_{\widetilde
s}^{\ssup2}
\bigr\}
\nonumber\\[-8pt]\\[-8pt]\nonumber
&&{}\times \bigl[\P_{1,2} \bigl(G_{N,s,\delta}^{\ssup1}\cap
G_{N,s,\delta
}^{\ssup2}\cap G_{N,\widetilde s,\delta}^{\ssup1}\cap
G_{N,\widetilde s,\delta}^{\ssup2} \bigr)
\\
&&{} -\P_{1,2} \bigl(G_{N,s,\delta
}^{\ssup1}
\cap G_{N,s,\delta}^{\ssup2} \bigr)\P_{1,2}
\bigl(G_{N,\widetilde s,\delta}^{\ssup1}\cap G_{N,\widetilde s,\delta
}^{\ssup2} \bigr)
\bigr].\nonumber
\end{eqnarray}

We now show, for any $s\neq\widetilde s$, that the limit superior of
the term in the last line is not positive. This finishes the proof by
Lebesgue's theorem.\vspace*{1pt}

We abbreviate $x=X^{\ssup1}_s$ and $\widetilde x=X^{\ssup
1}_{\widetilde s}$ and $y=X^{\ssup2}_s$ and $\widetilde y=X^{\ssup
2}_{\widetilde s}$. Without loss of generality, we assume that
$s<\widetilde s$, $x\neq y$ and $\widetilde x\neq\widetilde y$.
Furthermore we also may and will assume that $x\mathop
{\longleftrightarrow}\limits^{\geq}_{s} y$ and $\widetilde x\mathop
{\longleftrightarrow}\limits^{\geq}_{\widetilde s} \widetilde y$. Without
loss of generality, all the four terms $f_s(x),f_s(y),f_{\widetilde
s}(\widetilde x)$ and $f_{\widetilde s}(\widetilde y)$ are larger than
$\lambda_{\rm c}(R)$.
Let, as in the proof of Lemma~\ref{lem-Penrose}, $\Ccal_{x,\delta
}^{\ssup{s,N}}$ denote the biggest component of the union of the
$RN^{-1/d}$-balls around $X_s^{\ssup1}, X_s^{\ssup2},\dots
,X_s^{\ssup N}$ within $x+(-\delta,\delta)^d$, analogously for
$y,\widetilde s, \widetilde x$ and $\widetilde y$.

We recall from the proof of Lemma~\ref{lem-Penrose} that the
probability of the symmetric difference between $G_{N,t,\delta}^{\ssup
i}$ and the event $\{X_t^{\ssup i}\in\Ccal_{X_t^{\ssup i},\delta
}^{\ssup{t,N}}\}$, $i=1,2$ and $t=s,\widetilde s$, tends to $0$
as $N$ goes to infinity, followed by $\delta\downarrow0$. This
reduces the problem to showing that
%
\begin{eqnarray}
\label{eq-goal}
&& \limsup_{\delta\downarrow0} \limsup_{N\to\infty}
\bigl[\P _{1,2} \bigl(x\in\Ccal_{x,\delta}^\ssup{s,N}, y
\in\Ccal_{y,\delta
}^\ssup{s,N}, \widetilde x\in
\Ccal_{\widetilde{x}, \delta}^\ssup{\widetilde{s},N}, \widetilde y\in
\Ccal_{\widetilde y,\delta
}^{\ssup{\widetilde{s},N}} \bigr)
\nonumber\\[-8pt]\\[-8pt]\nonumber
&&\quad{} -\P_{1,2} \bigl(G_{N,s,\delta}^{\ssup1}\cap
G_{N,s,\delta}^{\ssup
2} \bigr)\P_{1,2} \bigl(G_{N,\widetilde s,\delta}^{\ssup1}
\cap G_{N,\widetilde s,\delta}^{\ssup2} \bigr) \bigr] \le0.
\end{eqnarray}

We pick $\delta>0$ smaller than $\frac{1}{3} \min\{\llvert  x-y\rrvert,\llvert  \widetilde
x-\widetilde y\rrvert  \}$. Let us give some heuristic explanation of the
following argument. To\vspace*{1pt} get \eqref{eq-goal}, we only have to prove
that, with probability tending to 1 as $N\to\infty$, the partial
clusters $\Ccal_{x,\delta}^{\ssup{s,N}}\cup\,\Ccal_{y,\delta
}^{\ssup{s,N}}$, and $\Ccal_{\widetilde x,\delta}^{\ssup
{\widetilde{s},N}}\cup\,\Ccal_{\widetilde y,\delta}^{\ssup{\widetilde{s},N}}$, depend only on two disjoint sub-collections of $X^{\ssup
3},\dots,X^{\ssup N}$ or at least on sub-collections with a small
overlap. What we mean precisely here is that the density of the walkers
in $\Ccal_{\widetilde x,\delta}^{\ssup{\widetilde{s},N}}\cup\,\Ccal
_{\widetilde y,\delta}^{\ssup{\widetilde{s},N}}$ is roughly the same
if we remove those points that were in $\Ccal_{x,\delta}^{\ssup{s,N}}\cup\,\Ccal_{y,\delta}^{\ssup{s,N}}$. Therefore,\vspace*{1pt} we need
Assumption~\ref{Ass-move}(ii) to describe the position of the walkers
at time $s$, given their position at
time $\widetilde{s}$. In more technical terms, it says the following. By
$\Bcal(D)$ we denote the Borel $\sigma$-field on $D$. Let a version
of the conditional distribution of $X_s$ given $X_{\widetilde s}=y$ be
given, that is, a Markov kernel $K_{s,\widetilde s}\colon D\times\Bcal
(D)\to\Bcal(D) $ such that, almost surely, $\P(X_s\in A\mid
X_{\widetilde s}=y)=K_{s,\widetilde s}(y,A)$ for any $A\in\Bcal(D)$.
Then we require that $K_{s,\widetilde s}(y,\{x\})=0$ for any $x\in D$.
Indeed, this assumption implies that, for any $y\in D$,
%
\begin{equation}
\label{Kerneldelta} \lim_{\delta\downarrow0}\P\bigl(X_s\in B(x,
\delta)\mid X_{\widetilde s}=y\bigr)= \lim_{\delta\downarrow0}K_{s,\widetilde s}
\bigl(y,B(x,\delta)\bigr) =K_{s,\widetilde s}\bigl(y,\{x\}\bigr)=0.
\end{equation}
Since the probability on the left-hand side is continuous in $y$ and
monotonous in $\delta$, the convergence is even uniform in $y\in D$,
according to Dini's theorem. Hence, we can multiply this term with
$f_{\widetilde s}(y)$, integrate over $y\in D$ and interchange this
integration with the limit $\delta\downarrow0$. Now we can see
heuristically the statement as follows. According to a large-$N$
ergodic theorem, there are only of order $N\delta^{2d}$ walkers that
are at time $s$ in $B(x,\delta)$ and at time $\widetilde s$ in
$B(\widetilde x,\delta)$, analogously with $y$ and $\widetilde y$.
Hence, among all the $\asymp N\delta^d$ walkers present in
$B(\widetilde x,\delta)$ at time $\widetilde s$, those ones who were
in $B(x,\delta)$ at time $s$ are negligible for small $\delta$. This
implies the claimed asymptotic independence.

Let us turn to the proof. We need to introduce a bit of notation. For
$A\subset\{1,\dots,N\}$, we write $\Ccal_{x,\delta}^\ssup{s,A}$
for the largest cluster in the $\delta$-box around $x$ that is built
out of all the $X^{\ssup i}_s$ with $i\in A$ only. We put
\[
A_s^{\ssup N}=\bigl\{i\in\{1,\dots,N\}\colon
X_s^{\ssup i}\notin B(x,\delta)\cup B(y,\delta)\bigr\}.
\]
Now we use the triangle inequality to bound
%
\begin{eqnarray}\label{thinning}
&& \P_{1,2} \bigl(x\in
\Ccal_{x,\delta}^\ssup{s,N}, y\in\Ccal _{y,\delta}^\ssup{s,N},\widetilde x\in\Ccal_{\widetilde x,\delta
}^\ssup{\widetilde{s},N},
\widetilde y\in\Ccal_{\widetilde y,\delta
}^\ssup{\widetilde{s},N} \bigr)\nonumber
\\
&&\quad \leq\P_{1,2} \bigl(x\in\Ccal_{x,\delta}^\ssup{s,N}, y
\in\Ccal _{y,\delta}^\ssup{s,N},\widetilde x\in
\Ccal_{\widetilde{x},\delta
}^{\ssup{\widetilde{s},A_s^{\ssup N}}}, \widetilde{y}\in\Ccal _{\widetilde{y},\delta}^{\ssup{\widetilde{s},A_s^{\ssup N}}}
\bigr)
\\
&&\qquad{}+\P_{1,2} \bigl(\widetilde x\in\Ccal_{\widetilde{x},\delta
}^{\ssup{\widetilde{s},N}}
\setminus\Ccal_{\widetilde{x},\delta
}^\ssup{\widetilde{s},A_s^{\ssup{N}}}
\bigr) +\P_{1,2} \bigl(\widetilde y\in\Ccal_{\widetilde{y},\delta}^{\ssup{\widetilde{s},N}}
\setminus\Ccal_{\widetilde{y},\delta}^\ssup{\widetilde{s},A_s^{\ssup{N}}}
\bigr).\nonumber
\end{eqnarray}
Since\vspace*{1pt} $\Ccal_{x,\delta}^\ssup{s,N}$ and $\Ccal_{y,\delta}^\ssup{s,N}$ depend only on the $X_s^{\ssup i}$ with $i$ in the complement
of $A_s^{\ssup N}$, the first two events in the first term on the
right-hand side are independent from the last two events. Lemma~\ref
{lem-Penrose}(ii) and the continuity of $\overline\theta(\cdot, R)$
imply that the probability of the intersection of the first two events
converges toward $\overline\theta(f_s(x),R)\overline\theta(f_s(y),R)$.
Note that the particles that the point processes $\Ccal_{\widetilde
{x},\delta}^\ssup{\widetilde{s},A_s^{\ssup{N}}}$ and $\Ccal
_{\widetilde{y},\delta}^\ssup{\widetilde{s},A_s^{\ssup{N}}}$ puts
are given by trajectories that do not visit any of the two balls
$B(x,\delta)$ and $B(y,\delta)$ at time $s$; more precisely, they are
picked according to the density
%
\begin{eqnarray}\label{fdensity}
f_{\widetilde s}^{\ssup{s,\delta}}(z)&=&\P\bigl(X_s
\notin B(x,\delta)\cup B(y,\delta),X_{\widetilde s}\in\d z\bigr)/\d
z
\nonumber\\[-8pt]\\[-8pt]\nonumber
&=& K_{s, \widetilde s} \bigl(z,\bigl(B(x,\delta)\cup B(y,\delta)
\bigr)^{\rm c} \bigr) f_{\widetilde s}(z).
\end{eqnarray}
Hence, the probability of the intersection of the last two events
converges toward
\[
\overline\theta\bigl(f_{\widetilde s}^{\ssup{s,\delta}}(\widetilde x),R\bigr)
\overline\theta\bigl(f_{\widetilde s}^{\ssup{s,\delta}}(\widetilde y),R\bigr).
\]

A glance at \eqref{fdensity} shows that $f_{\widetilde s}^{\ssup{s,\delta}}(z)$ converges, as $\delta\downarrow0$, for any $z\in D$,
toward $\P(X_s\neq x,X_s\neq y,X_{\widetilde s}\in\d z)/\d z$, which
is, by Assumption~\ref{Ass-move}(i) (or also by (ii)), equal to
$f_{\widetilde s}(z)$.
Since $f_{\widetilde s}(\widetilde x)$ and $f_{\widetilde s}(\widetilde
y)$ are larger than the critical value, we may use continuity of
$\overline\theta$.

All together, we have that the first term of the right-hand side of
\eqref{thinning} converges, as $N\to\infty$ followed by $\delta
\downarrow0$, toward
%
\begin{equation}
\label{eq-indeppairs} \overline\theta\bigl(f_s(x),R\bigr)\overline\theta
\bigl(f_s(y),R\bigr)\overline\theta \bigl(f_{\widetilde s}(
\widetilde x),R\bigr)\overline\theta\bigl(f_{\widetilde
s}(\widetilde{y}),R
\bigr).
\end{equation}

Furthermore, Assumption~\ref{Ass-move}(ii) also implies that
%
\begin{equation}
\label{eq-missing} \limsup_{N\to\infty}\P_{1,2} \bigl(
\widetilde x\in\Ccal _{\widetilde{x},\delta}^{\ssup{\widetilde{s},N}}\setminus\Ccal
_{\widetilde{x},\delta}^\ssup{\widetilde{s},A_s^{\ssup{N}}}
\bigr)
\end{equation}
vanishes as $\delta\downarrow0$. Indeed, we know that $\Ccal
_{\widetilde{x},\delta}^\ssup{\widetilde{s},A_s^{\ssup{N}}}\subset
\Ccal_{\widetilde{x},\delta}^{\ssup{\widetilde{s},N}}$, therefore
the above limit superior is equal to
$\overline\theta(f_{\widetilde{s}}(\widetilde{x}))-\overline\theta
(f_{\widetilde s}^{\ssup{s,\delta}}(\widetilde{x}))$. Hence, the
convergence of $f_{\widetilde s}^{\ssup{s,\delta}}$ and the
continuity of $\overline\theta$ give the result. We proceed
analogously for the last term in \eqref{thinning} and get that the
limit superior as $N\to\infty$ and $\delta\downarrow0$ of the
left-hand side of \eqref{thinning} is not larger than the expression
in \eqref{eq-indeppairs}. Now use Lemma~\ref{lem-Penrose}(ii) for the
second term in \eqref{eq-goal} to see that from this the desired
assertion follows.
\end{pf}

\begin{pf*}{Proof of Theorem~\ref{thm-connectiontime}}
First note that both assertions of \eqref{taulimit} easily follow from
Proposition~\ref{prop-exptau}, in conjunction with Lemma~\ref
{lem-prob-conv}, provided that $R$ is chosen such that
%
\begin{equation}
\label{X1notequalX2} \int_0^T \d s \1\bigl\{
f_s\bigl(X^{\ssup i}_s\bigr)=
\lambda_{\rm c}(R)\bigr\} =0\qquad \mbox{for }i=1,2.
\end{equation}
Furthermore, note that, almost surely, \eqref{X1notequalX2} holds for
almost all $R$. Indeed, this follows from
\begin{eqnarray*}
&& \E \biggl(\int_0^{\infty} \d R \int_0^T\d s \1 \bigl\{f_s
\bigl(X_s^{\ssup
i}\bigr)=\lambda_{\rm c}(R)\bigr\}
\biggr)
\\
&&\quad = \int_0^T \d s \int_D
\,\d x f_s(x) \int_0^{\infty} \d R \1
\bigl\{ f_s(x)=R^{-d} \lambda_{\rm c}(1)\bigr\} =0.
\end{eqnarray*}
Hence, for a given (random) exceptional $R$, we pick sequences
$(R_k)_{k\in\N}$ and $(R_k')_{k\in\N}$ such that $R_k\downarrow R$
and $R_k'\uparrow R$ and $R_k$ and $R_k'$ satisfy \eqref{X1notequalX2}
for any $k$ in place of $R$. Since $\tau_T^{\ssup N}$ is an increasing
function of $R$, we may estimate it from above and below by replacing
$R$ with $R_k$ and $R_k'$, respectively,\vspace*{1pt} and applying Proposition~\ref
{prop-exptau} and Lemma~\ref{lem-prob-conv} with these. This yields
\eqref{taulimit} with $\tau_T^{\ssup\geq}$ and $\tau_T^{\ssup>}$
replaced by their versions for $R$ replaced with $R_k$ and with $R_k'$,
respectively.

The only thing that we need to do is to show the right-upper
semicontinuity of the map $R\mapsto\tau_T^{\ssup\geq}$ and the
left-lower semicontinuity of the map $R\mapsto\tau_T^{\ssup>}$. To
show these, note that $\overline\theta^{\ssup\geq}(\cdot
,R)=\overline\theta(R^d\cdot +,1)$ is right-continuous and
$\overline\theta^{\ssup>}(\cdot,R)=\overline\theta(R^d\cdot
-,1)$ is left-continuous. Furthermore, for any $x,y\in D$ and any $s\in
[0,T]$, the map $R\mapsto\1\{x\mathop{\longleftrightarrow}\limits^{\geq
}_{s}y\}$ is\vspace*{-2pt} right-upper semicontinuous, and the map $R\mapsto\1\{
x\mathop{\longleftrightarrow}\limits^{>}_{s}y\}$ is left-lower
semicontinuous. The latter assertion is quite easy to see; let us show
the former. Assume that, for all $\eps>0$, $x$ and $y$ are connected
through the set $\{f_s\geq\lambda_{\rm c}(R+\eps)\}$. Recall that
$\lambda_{\rm c}(R)=R^{-d}\lambda_{\rm c}(1)$ is decreasing in $R$.
If $x$ and $y$ were not connected through the set $\{f_s\geq\lambda
_{\rm c}(R)\}$, then they would lie in different components
of this set. By compactness, these components have a positive distance
to each other. Hence, there is a hyperplane in $D$ through the
complement of $\{f_s\geq\lambda_{\rm c}(R)\}$ that separates these
two components. Since this hyperplane is compact, $f_s$ assumes a
maximum on it, which is strictly smaller than $\lambda_{\rm c}(R)$.
Hence, every curve from $x$ to $y$ must cross this hyperplane, that is,
must pass a point with an $f_s$-value bounded away from $\lambda_{\rm
c}(R)$. This means that, for some sufficiently small $\eps>0$, $x$ and
$y$ are not connected through $\{f_s\geq\lambda_{\rm c}(R+\eps)\}$.
Hence, $\limsup_{\eps\downarrow0}\1\{x\mathop{\longleftrightarrow
}\limits^{\geq, R+\eps}_{s}y\}\leq\1\{x\mathop{\longleftrightarrow}\limits^{\geq
, R}_{s}y\}$, where we wrote $\mathop{\longleftrightarrow}\limits^{\geq,
R}_{s}$ for connectedness through the set $\{f_s\geq\lambda_{\rm
c}(R)\}$. Using Lebesgue's theorem shows the claimed continuity
properties of $\tau_T^{\ssup\geq}$ and $\tau_T^{\ssup>}$ in $R$
and finishes the proof of Theorem~\ref{thm-connectiontime}.
\end{pf*}

\section{Long-time investigations for the random waypoint model}\label
{sec-prooflongtime}

  In this section, we prove Lemma~\ref{lem-ergodic} and
Theorem~\ref{thm-longtimebound}, that is, we restrict ourselves to the
random waypoint model (RWP) introduced in Section~\ref{sec-RWPResults}
and study the long time behaviour of the limiting connection time both
in terms of an ergodic theorem and a large-deviations result.
First, we prove in Section~\ref{invariant-distrib} the convergence of
the RWP to its invariant distribution. The proof of Lemma~\ref
{lem-ergodic} is based on a certain discrete-time Markov chain, whose
ergodic and mixing properties are derived in Section~\ref{sec-Mixing}.
The proof then follows in Section~\ref{sec-ergodicproofII}. Finally,
we prove Theorem~\ref{thm-longtimebound} in Section~\ref{sec-proofThDeviat}.

\subsection{Recurrence and ergodicity of the RWP}
\label{invariant-distrib}

  Since we want to study long-time properties of the connection
time, we will need recurrence and ergodic properties of the RWP, which
we provide in this section. For the special case of $\Wcal$ being the
uniform distribution on $D$, most of our results in this section are
already contained in \cite{BoudecVojnovic06}, but our Proposition~\ref
{Harris-rec} below also contains a statement on convergence in total
variation, which will be important in Lemma~\ref{lem-ZHarris} below.
For the reader's convenience, we provide all necessary proofs; they are
independent of \cite{BoudecVojnovic06}, but use different variants of
the Markov renewal theorem available in the literature.

The trajectory is divided into {\it trips}, by which we mean the parts
from leaving a waypoint to arriving at the next one. $\P^{\ssup0}$
and $\E^{\ssup0}$ denote probability and expectation if the process
starts at time 0 at the beginning of a trip at the zeroth of the
waypoints, that is, if the initial waypoint $W_0$ has distribution
$\Wcal$.

In \cite{BoudecVojnovic06}, Theorem~6, another variant of $Y$ is
considered, and it is argued that process possesses a unique invariant
distribution. Projecting on our first coordinate, the location of the
walker, the distribution of $X$ in equilibrium is given by the formula
%
\begin{equation}
\label{invmeasX} \mu_*(\d x)=\frac{1}{Z}\int_0^1
\d s \E^{\ssup0} \biggl(\frac
{V_1}{\llvert  W_1-W_0\rrvert  };W_0+s(W_1-W_0)
\in\d x \biggr),
\end{equation}
where $Z$ is a normalisation. It turns out below that this formula
persists also for a general waypoint measure. In particular $\mu_*$
has a continuous density. We refer in particular to \cite{Boudec04}
for a general methodology to describe this measure. See \cite{BettstetterWagner02}, Section~5,
and \cite{HyytiaLassilaVirtamo06}, Sections~III and~IV, for explicit formulas, approximations and
simulations for special cases of domains $D$ and waypoint measures
$\Wcal$, like uniform distributions on rectangles and balls.

For the sake of illustration, we give an explicit value in $d=2$ in the
simplest case where the domain is the unit disk, the waypoint measure
$\Wcal$ is the uniform measure on it and the velocity is chosen to be
constant. In this case, the density of the waypoint location in the
invariant distribution is given by\vspace*{-1pt}
\[
f_*(x)=\frac{45}{64 \pi}\bigl(1-\llvert x\rrvert ^2\bigr)\int
_0^\pi\sqrt{1-\llvert x\rrvert ^2
\cos ^2(\varphi)} \,\d\varphi,\qquad x\in B(0,1). %
\]
An approximation with a mean square error $\leq0.0065$ and an absolute
error $\leq0.067$ is given by $f_*(x)=\frac{2}{\pi} (1-\llvert  x\rrvert  ^2)$; see
\cite{QZ07} and \cite{BettstetterWagner02}, equation~(18).

In the following, we give detailed proofs for ergodic properties of the
RWP, based on the Markov renewal theorem in the form provided by \cite
{Kesten74}. Alternative proofs could be based on the form given in
\cite{BoudecVojnovic06}, Theorem~6.

We first show that the sequence of the trips is positive Harris
recurrent. More precisely, we consider the sequence $\mathcal
T=(\mathcal T_n)_{n\in\N}=(W_{n-1},W_{n},V_{n})_{n\in\N}$
in $\Dcal$.
Since $(W_n)_{n\in\N_0}$ and $(V_n)_{n\in\N}$ are independent
i.i.d.~sequences, $\mathcal T$ is obviously a Markov chain.
Furthermore, it is also easy to see that $\mathcal T$ is positive
Harris recurrent, since it satisfies
%
\begin{equation}
\label{THarris} \P_y(\mathcal T_n\in A)=\Wcal\otimes
\Wcal\otimes\Vcal(A),\qquad n\geq2,y\in\Dcal,A\subset\Dcal\mbox{ mb.},
\end{equation}
where we wrote $\P_y$ for the probability measure under which the
walker starts from $Y(0)=y$. We use this to prove the convergence of
$Y_t$ introduced in \eqref{Ydef}. The proof goes in two step. The
first one (see Lemma~\ref{lem-Harris}) applies the Markov renewal
theorem using the fact that $Y_t$ is a time change of $\mathcal T$ and
gives a good understanding and a description of the limit law (in
particular it states the existence of an invariant distribution with
finite mass). However, as we will see, this approach only gives weak
convergence. In a second step, we use Harris recurrence (see
Proposition~\ref{Harris-rec}) to obtain convergence in total
variation. Of course it is then easy to check that the convergence has
to be toward the same limit. By $\P_\alpha$, we denote the
probability measure under which the process $(Y_t)_{t\in[0,\infty)}$
starts from the distribution $\alpha$.

\begin{lemma}\label{lem-Harris}
For any bounded continuous function $g\colon\Dcal\times\R^+\to\R
^+$, and for any $y\in\Dcal$,
%
\begin{equation}
\lim_{t\to\infty} \E_y\bigl[g(\mathcal{T}_{N(t)},T_{N(t)}-t
)\bigr]=\frac
{1}{\E[U_1]}\int_{\Dcal}\P_{\Wcal\otimes\Wcal\otimes\Vcal
}[
\mathcal{T}_1\in\d z, U_1 \in\d\lambda]\int
_0^\lambda g(z,s)\, \d s.
\end{equation}
\end{lemma}

\begin{pf}
We apply \cite{Kesten74}, Theorem~1, which immediately implies the
assertion, noting that the measure $\psi$ in \cite{Kesten74} is
indeed equal to $\Wcal\otimes\Wcal\otimes\Vcal$ by \cite{Kesten74}, Lemma~2, that is, we only have to check the validity of Conditions
I.1--I.4 of \cite{Kesten74}.

Conditions I.1 and I.2 are trivial here, while Condition I.3 is the
usual non-lattice assumption. It states that there is a non-lattice
sequence $(\zeta_\nu)_{\nu\in\N}$ in $\R$ such that, for each
$\nu\in\N$ and $\delta>0$, there exists some $y\in\Dcal$, such
that, for every $\epsilon>0$, there exists a measurable set $A$ with
positive $\Wcal\otimes\Wcal\otimes\Vcal$-measure, integers $m_1$,
$m_2$ and $\tau\in\R$ such that, for $x\in A$,\vspace*{-2pt}
%
\begin{eqnarray}\label{Cond3}
\P_x\bigl[d(\mathcal{T}_{m_1},y)<\epsilon,
\llvert T_{m_1}-\tau\rrvert \le\delta \bigr]&>&0\quad\mbox{and}
\nonumber\\[-8pt]\\[-8pt]\nonumber
\P_x\bigl[d(\mathcal{T}_{m_2},y)<\epsilon, \llvert
T_{m_2}-\tau-\zeta_\nu\rrvert \le\delta\bigr]&>&0,
\end{eqnarray}
$d$ being the usual Euclidean distance on $\Dcal$.\vadjust{\goodbreak}

We will prove this assumption with an arbitrary $y=(w_0,w_1,v_1)$
inside the support of $\Wcal\,\otimes\,\Wcal\,\otimes\,\Vcal$, not
depending on $\nu$ nor on $\delta$, and with $A=\{x\in\Dcal\colon
d(x,y)<\epsilon\}$, where we assumed without loss of generality that
$2\epsilon v_{-}^{-1}+\diam(D)\epsilon v_-^{-2}<\delta/3$.
Furthermore, we put $\tau:=\llvert  w_1-w_0\rrvert  /v_1$ and pick any non-lattice
sequence $(\zeta_\nu)_{\nu\in\N}$ inside the support of $\tau
+\llvert  w_0-w_1\rrvert  /V_1$. Furthermore, put $m_1=1$ and $m_2=3$. By continuity of
the densities of $\Wcal$ and $\Vcal$, the $\Wcal\,\otimes\,\Wcal\,\otimes\,\Vcal$-measure of $A$ is positive. Putting
$x=(w_0',w_1',v_1')\in A$ and denoting by $T_1(x)=\llvert  w_1'-w_0'\rrvert  /v_1'$
the (deterministic) value of $T_1$ starting from $x$, we see that
%
\begin{eqnarray}\label{Rxcond}
\bigl\llvert T_1(x)-\tau\bigr\rrvert &\leq&
\frac{\llvert  w_1'-w_0'-(w_1-w_0)\rrvert  }{v_1'}+\llvert w_1-w_0\rrvert \biggl\llvert
\frac{1}{v_1}-\frac{1}{v_1'}\biggr\rrvert
\nonumber\\[-8pt]\\[-8pt]\nonumber
&\leq& \frac{2\epsilon}{v_-}+
\frac{\diam(D)\epsilon}{v_-^2}<\frac{\delta}3.
\end{eqnarray}
Noting that $\Tcal_1=x$ with $\P_x$-probability one, we see that the
first part of \eqref{Cond3} is satisfied; the probability is even
equal to one.

Now we turn to the proof of the second. Keep $x\in A$ fixed. Recall
that $T_n=U_0+U_1+\cdots+U_{n-1}$ and that $U_n=\llvert  W_{n+1}-W_{n}\rrvert  /V_n$
for any $n$. Note that, under $\P_x$, $\mathcal T_3$ has distribution
$\Wcal\,\otimes\,\Wcal\,\otimes\,\Vcal$, and therefore $\P_x(d(\mathcal
{T}_{3},y)<\epsilon)=\Wcal\,\otimes\,\Wcal\,\otimes\,\Vcal(A)>0$. On the
event $\{d(\mathcal{T}_{3},y)<\epsilon\}$, with $\P_x$-probability
one, \eqref{Rxcond} shows that $\llvert  U_0-\tau\rrvert  <\delta/3$, and a the same
calculation with $x$ replaced by $\mathcal T_3$ shows that $\llvert  U_2-\tau
\rrvert  <\delta/3$. By our choice of $\zeta_\nu$ and by continuity of the
densities of $\Wcal$ and $\Vcal$, we easily see that the event $\{
\llvert  U_1+\tau-\zeta_\nu\rrvert  \leq\delta/3 \}$ has positive $\P
_x$-probability on $\{d(\mathcal{T}_{3},y)<\epsilon\}$, since
\begin{eqnarray*}
\llvert U_1+\tau-\zeta_\nu\rrvert
&\leq& \frac{\llvert  W_2-W_1-(w_0-w_1)\rrvert  }{v_-}+\biggl\llvert \frac{\llvert  w_0-w_1\rrvert  }{V_2}-(\zeta_\nu-
\tau)\biggl\llvert
\\
&\leq&\frac{2\epsilon
}{v_-}+\biggr\rrvert \frac{\llvert  w_0-w_1\rrvert  }{V_2}-(
\zeta_\nu-\tau)\biggr\rrvert,
\end{eqnarray*}
and the probability (with respect to $V_2$) to have the last term
smaller than $\diam(D)\epsilon v_-^{-2}$ is positive. Since
\[
\llvert T_{3}-\tau-\zeta_\nu\rrvert =\llvert
U_0+U_1+U_2-\tau-\zeta_\nu\rrvert
\leq\llvert U_0-\tau \rrvert +\llvert U_1+\tau-
\zeta_\nu\rrvert +\llvert U_2-\tau\rrvert, %
\]
we now see that also the last condition in \eqref{Cond3} is satisfied.

Condition I.4 states that, for any $x\in\Dcal$, $\delta>0$, there
exists $r_0(x,\delta)>0$ such that for any measurable function
$f\colon\Dcal^{\N}\times\R^{\N_0}\to\R$, and for all $y$ with
$d(y,x)<r_0(x,\delta)$,
%
\begin{eqnarray}\label{assumptioni4}
&&\E_x \bigl[f \bigl((\mathcal{ T}_i)_{i\in\N},(U_i)_{i\in\N_0}
\bigr) \bigr]
\nonumber\\[-8pt]\\[-8pt]\nonumber
&&\quad
\le\E_y \Bigl[ \lim_{n\to\infty} \sup\bigl\{f
\bigl((t_i)_{i\in\N}, (u_i)_{i\in\N_0}
\bigr)\colon d(t_i,\mathcal {T}_i)+\llvert
u_i-U_i\rrvert <\delta\mbox{ for }i\le n\bigr\}
\Bigr]+\delta\sup\llvert f\rrvert.
\end{eqnarray}

This assumption is in general difficult to prove, but here things are
simple, as $\mathcal{T}_i$ and $U_i$ are independent of the starting
point for $i\ge3$. We can do the following coupling: write
$x=(w_0^\ssup{x},w_1^\ssup{x},v_1^\ssup{x})$ and $y=(w_0^\ssup
{y},w_1^\ssup{y},v_1^\ssup{y})$. We draw a sequence of
i.i.d.~waypoints and speeds $(W_i,V_i)_{i\ge2}$ according to $\Wcal\,\otimes\,\Vcal$. Define, for $z\in\{x,y\}$,
%
\begin{equation}
{W}_0^\ssup{z}=w_0^\ssup{z},
\qquad{W}_1^\ssup{z}=w_1^\ssup{z},
\qquad{V}_1^\ssup{z}=v_1^\ssup{z},\qquad
\bigl({W}_i^\ssup {z},{V}_i^\ssup{z}
\bigr)_{i\ge2}=(W_i,V_i)_{i\ge2},
\end{equation}
and put $\mathcal{T}^{\ssup{z}}_i=(W^{\ssup{z}}_{i-1},W^{\ssup
{z}}_i,V^{\ssup{z}}_i)$. It is then clear that $ (\mathcal{T}^{\ssup
{z}}_i)_{i\in\N}$ is a realisation of $(\mathcal{T}_i)_{i\in\N}$
under $\P_z$ and that for any $i\ge3$,
$\mathcal{T}^{\ssup{x}}_i=\mathcal{T}^{\ssup{y}}_i$. We saw in the
verification of Condition~I.3 that, if $d(x,y)<r$, then with obvious notation,
\[
d\bigl(\mathcal{T}_i^\ssup{x},\mathcal{T}_i^\ssup{y}
\bigr)<r, \qquad d\bigl(U_i^\ssup{x},U_i^\ssup{y}
\bigr)<r \biggl(\frac{2}{v_-}+\frac{\diam
(D)}{v_-^2} \biggr). %
\]
Taking $r_0(\delta)$ such that both right-hand sides are $<\delta$,
immediately gives Condition I.4.
\end{pf}

Using \eqref{Xtformula}, we easily derive the above mentioned weak
convergence of $X_t$ toward $\mu_*$ identified in \eqref{invmeasX},
as $X_t$ may be written as an explicit continuous function of $\mathcal
{T}_{N(t)} $ and $ T_{N(t) }-t$. We now give a refined result, using
the notion of Harris recurrence for continuous-time Markov chains.
First, note that the process
\[
\mathcal{Y}= (\mathcal{Y}_t )_{t\in[0,\infty)}= \biggl(
\mathcal{T}_{N(t)},\frac{T_{N(t)}-t}{U_{N(t)-1}} \biggr)_{t\in
[0,\infty)} %
\]
is a continuous-time Markov chain on $\Dcal\times[0,1]$ with
right-continuous paths. The second component of $\mathcal{Y}$ runs
from 0 to 1 with linear speed between the arrival times at the
waypoints. It is also easy to express $Y_t$ as a continuous functional
of $\mathcal{Y}_t$.

\begin{prop}\label{Harris-rec} $(\mathcal{Y}_t)_{t\in[0,\infty)}$
is a strongly aperiodic Harris recurrent chain, and its distribution
converges in total variation toward the unique invariant distribution.
As a consequence, the convergence in Lemma~\ref{lem-Harris} is true
for any measurable bounded function $g$. Furthermore, an ergodic
theorem holds for $(\mathcal{Y}_t)_{t\in[0,\infty)}$.
\end{prop}

\begin{pf}
We use the characterisation of Harris recurrence given in
\cite{KaspiMandelbaum94}, Theorem~1, with the measure $\nu$ given by
$\Wcal\otimes\Wcal\otimes\Vcal\otimes\lambda$, where $\lambda$
is the Lebesgue measure on $[0,1]$. It is easy to see that any set $A$
with positive $\nu$-measure will be hit by the process $(\mathcal
{Y}_t)_{t\in[0,\infty)}$. Indeed, without loss of generality, we can
assume that $A$ is a product set. By independence it will certainly
happen that one of the $\mathcal{T}_n$ will fall into the $\Dcal
$-component of $A$. Then as $\frac{T_{N(t)}-t}{U_{N(t)-1}}$ visits all\vspace*{1pt}
of $[0,1]$ between two waypoints, it follows that also $A$ will be hit
by $\mathcal{Y}$, implying Harris recurrence.

This implies in particular the existence of a unique (up to
multiplicative constants) invariant measure. It is not difficult to
check that this measure has to be the one appearing in Lemma~\ref
{lem-Harris}, up to the normalisation. In particular, it has finite
total mass. As a consequence, $\mathcal{Y}$ is strongly Harris
recurrent. We also have that this process has spread-out cycles, in the
sense of \cite{A03}, page~202. In fact, the hitting times of any set
under any starting point are spread out. Indeed, the first hitting
times might be deterministic (if the initial condition implies that the
set is hit during the first travel of the walker), but then one can
easily check that, due to the existence of a density for the speed, the
hitting times also have a continuous density. Therefore, using \cite{A03}, Proposition VII.3.8, this implies convergence in total variation
of $\mathcal Y_t$ toward its invariant distribution. The ergodic
theorem can be found in \cite{A03}, Proposition~VII.3.7.
\end{pf}

Note that, at this point, it would be possible to use the above result
to get a simple proof of Lemma~\ref{lem-ergodic}. However, we would
like to present a different proof, as we need to introduce the
important discrete-time Markov chain $(Z_j)_{j\in\N_0}$, that will be
useful for the sequel. This proof can be found in Section~\ref
{sec-ergodicproofII}.

\subsection{Recurrence and mixing properties of $Z$}\label{sec-Mixing}

  In this section, we introduce an important tool for our
proofs of Lemma~\ref{lem-ergodic} and Theorem~\ref
{thm-longtimebound}, a discrete-time Markov chain $Z=(Z_j)_j$ that
registers the locations, waypoints and velocities of two independent
RWPs at all the times at which one of them arrives at a new waypoint.
In this section, we study recurrence and the mixing properties of this
chain, in Sections~\ref{sec-ergodicproofII} and~\ref
{sec-proofThDeviat} we will use it to derive the long-time average and
large-deviations properties of the connection time. For proving just
the former of the two results in Lemma~\ref{lem-ergodic}, some
straight-forward ergodic arguments would be also sufficient, however,
we will need the identification of the ergodic limit in terms of the
Markov chain $Z$ in order to prove the large-deviations result in
Theorem~\ref{thm-longtimebound}. We show that $(Z_k)_{k\in\N_0}$ is
a time-homogeneous, $\psi$-mixing and Harris ergodic Markov chain. It
is an object of independent
interest, as it may serve also for other long-time investigations of
the model, as well as for computer simulations.

The Markov chain $Z$ is defined as follows. We consider the times
$0\leq S_1<S_2<\cdots$ at which any of the two walkers arrives at his
waypoint. Formally, $S_0=0$ and
%
\begin{equation}
\label{Sdef} S_j=\inf\bigl\{t>S_{j-1}\colon
W_{N^{\ssup1}(t)}^{\ssup1}\neq W_{N^{\ssup1}(S_{j-1})}^{\ssup1}
\mbox{ or } W_{N^{\ssup
2}(t)}^{\ssup2}\neq W_{N^{\ssup2}(S_{j-1})}^{\ssup2}
\bigr\},\qquad j\in\N,
\end{equation}
where the superscripts $(1)$ and $(2)$ mark the two walkers. Put
%
\begin{eqnarray}
\label{Zdef} Z_j &=& \bigl(Y^{\ssup1}(S_j),Y^{\ssup2}(S_j)
\bigr)
\nonumber\\[-8pt]\\[-8pt]\nonumber
&=& \bigl( \bigl(X^{\ssup1}_{S_j},W_{N^{\ssup1}(S_j)}^{\ssup1},V_{N^{\ssup
1}(S_j)}^{\ssup1}
\bigr), \bigl(X^{\ssup2}_{S_j},W_{N^{\ssup
2}(S_j)}^{\ssup2},V_{N^{\ssup2}(S_j)}^{\ssup2}
\bigr) \bigr)\in\Dcal^2, \qquad j\in\N_0.
\end{eqnarray}

That is, $Z=(Z_j)_{j\in\N_0}$ is the trace-Markov chain of two
independent copies of the RWP, observed at the times at which any of
the two arrives at a waypoint; it is a time-change of $(Y^{\ssup
1},Y^{\ssup2})$. It is easy to see that $(Z_j)_j$ is a
time-homogeneous Markov chain on $\Dcal^2$. This chain does not
explicitly record the location of the random walker at any fixed time,
but the time that passes between the waypoint arrivals can be deduced
from the information contained in $Z$. Hence, it is well-suitable for
deducing asymptotic assertions for long time. First we derive a mixing
property, which will later be used for the large-deviations principle.

\begin{lemma}\label{lem-mixing}
The sequence $(Z_j)_j$ is $\psi$-mixing under any starting
distribution, that is,
\[
\lim_{k\to\infty} \sup_{A\in\mathcal F_0^0,B\in\mathcal
F_{k}^{\infty}} \biggl\llvert
\frac{\P(A\cap B)}{\P(A) \P(B)}-1\biggr\rrvert = 0, %
\]
where $\mathcal F_m^k:=\sigma(Z_m,\dots,Z_k)$.
\end{lemma}

\begin{pf}
Introduce the event
\begin{eqnarray*}
E_k&=&\bigl\{\exists l,m\in\{1,\dots,k-1\}: W_0^{\ssup
1}\neq W_{N^{\ssup1}(S_l)}^{\ssup1}\neq W_{N^{\ssup1}(S_k)}^{\ssup1}
\mbox{ and}\\
&&{} W_0^{\ssup2}\neq W_{N^{\ssup2}(S_m)}^{\ssup2}
\neq W_{N^{\ssup2}(S_k)}^{\ssup2}\bigr\}
\end{eqnarray*}
that both walkers choose at least two new waypoints by time $S_k$.
Then, conditional on $E_k$, any $A\in\mathcal F_0^0$ and $B\in
\mathcal F_{k}^{\infty}$ are independent. Indeed, on the event $E_k$,
$A$ depends on $X_0^{\ssup1}, W_1^{\ssup1}, V_1^{\ssup1},X_0^{\ssup
2}, W_1^{\ssup2}, V_1^{\ssup2}$ only, while $B$ depends\vspace*{1pt} only on the
variables $W_l^{\ssup1}$, $V_l^{\ssup1}$, $W_l^{\ssup2}$, $V_l^{\ssup2}$ for
some $l\geq3$ and on $X_{S_l}^{\ssup1}, X_{S_l}^{\ssup2}$ with
$l\geq2$; note that, for $i\in\{1,2\}$, $X_{S_l}^{\ssup i}$ is a
function of $W_{N^{\ssup i}(S_l)}^{\ssup i},W_{N^{\ssup
i}(S_l)-1}^{\ssup i}$ and $V_{N^{\ssup i}(S_l)}^{\ssup i}$ only, and
$N^{\ssup i}(S_l)\geq3$ on $E_k$. Using the independence of $A$ and
$B$ on $E_k$, a small calculation yields that
\[
\frac{\P(A\cap B)}{\P(A) \P(B)} =\frac{\P(E_k\llvert  A)\P(E_k\rrvert  B)}{\P(E_k)}+\P\bigl(E_k^{\rm c}
\mid A\cap B\bigr)\frac
{\P(A\cap B)}{\P(A) \P(B)}.
\]
Hence, the assertion follows from
%
\begin{equation}
\label{Ukconv} \lim_{k\to\infty}\sup_{A\in\mathcal F_0^0,B\in\mathcal
F_{k}^{\infty}}\P
\bigl(E_k^{\rm c}\mid A\cap B\bigr)=0.
\end{equation}
We show now that \eqref{Ukconv} holds. The event $E_k^{\rm c}$ splits
into the event that the first walker has chosen not more than one new
waypoint by time $S_k$, but the second has chosen at least $k-1$ new
waypoints, and the same event with first and second walker reversed.
Let us only look at the first of these two events. On this event, the
time $S_k$ is not larger than $2\diam(D)/v_-$, since a choice of a new
waypoint is done after $\diam(D)/v_-$ time units at the latest, since
all ways are no longer than $\diam(D)$ and all velocities are no less
than $v_-$. Since\vspace*{1pt} the time that passes between the second walker picks
his $(j-1)$st and the $j$th waypoint is $\llvert  W^{\ssup2}_{j}-W^{\ssup
2}_{j-1}\rrvert  /V^{\ssup2}_{j}$, we have\vspace*{2pt} that its sum over $j\in\{1,\dots
,k-1\}$ is not larger than $2\diam(D)/v_-$. Hence, on this event we have
\[
\sum_{j=1}^{k-1}\bigl\llvert
W^{\ssup2}_j-W^{\ssup2}_{j-1}\bigr\rrvert
\leq2\frac
{v_+}{v_-}\diam(D). %
\]
Leaving out the summands for $j=1$ and $j=k-1$, this remaining sum is
still upper bounded by the right-hand side, and it does not depend on
$Z_0$ nor on $Z_k,Z_{k+1},\dots.$ Hence, the probability for this sum
being smaller than the right-hand side is an upper bound for the half
of $\P(E_k^{\rm c}\mid A\cap B)$ that we are considering, and it does not
depend on $A$ nor on $B$. Since the right-hand side is constant and
since the waypoints are not deterministic, the probability for this
event tends to $0$ as $k\to\infty$. This shows that \eqref{Ukconv}
holds and ends the proof.
\end{pf}

The following lemma says that $Z$ is Harris recurrent, has a unique
invariant distribution and is non-lattice, which is summarised by
saying that it is Harris ergodic. In particular, it satisfies an
ergodic theorem, that is, for any bounded measurable function $f$, the
averages $\frac{1} N\sum_{i=1}^N f(Z_i)$ converge almost surely to the
integral of $f$ with respect to the invariant distribution.

\begin{lemma}\label{lem-ZHarris}
The chain $Z$ is Harris ergodic.
\end{lemma}

\begin{pf}
Harris recurrence of $Z$ is equivalent to the existence of a
non-trivial $\sigma$-finite measure $\varphi$ such that $Z$ is
$\varphi$-recurrent, see \cite{A03}, Corollary~VII.3.12. Therefore,
we have to show that there exists some $\sigma$-finite measure
$\varphi$ such that every measurable set $F\subset\mathcal D^2$ with
$\varphi(F)>0$ is recurrent.

We denote the invariant measure of the process $(Y^{\ssup1}_t)_{t\in
[0,\infty)}$ by $\gamma$. Define $\varphi= {\gamma}\otimes\Wcal
\otimes\Wcal\otimes{\Vcal}$, which is obviously $\sigma$-finite.
Let $F\subset\mathcal D^2$ be measurable with $\varphi(F)>0$. We are
going to show that the hitting time of $F$ is almost surely finite.
Note that $\varphi(F)>0$ implies, by Fubini's theorem, that, for some
$\epsilon>0$, the set $\widetilde F$ of all $T$ satisfying $\int\1
_F((Y,T)) \gamma(\d Y)>\epsilon$ has positive $\Wcal\otimes\Wcal
\otimes{\Vcal}$ measure.

First, consider the sequence $(n_k)_{k\in\N_0}$ of times at which the
second walker arrives at a waypoint, that is, $(S_{n_k})_{k\in\N
_0}=(T_k^{\ssup2})_{k\in\N_0}$. The first component of the process
$(Z_{n_k})_{k\in\N_0}$ is a RWP sampled at times which are given by
an independent renewal process, and the second component has the same
law as $({\mathcal T}_{k+1})_{k\in\N_0}$. According to \eqref
{THarris} and \cite{A03}, Corollary~VII.3.12, the second component is
$(\Wcal\otimes\Wcal\otimes{\Vcal})$-positive recurrent. In
particular there exists a subsequence $(\widetilde{n}_k)_k$ of
$(n_k)_k$ such that the second component of $Z_{\widetilde{n}_k}$
belongs to $\widetilde{F}$ for any $k\in\N_0$. Also $(S_{\widetilde
{n}_k})_{k\in\N_0}$ is a transient Markov renewal process,
independent of $Y^\ssup{1}$.

Now conditioning on the second component process, $Y^\ssup{2}$,
$(Y^\ssup{1}_{S_{\widetilde{n}_k}})_{k\in\N_0}$ is given by
sampling the, by Proposition~\ref{Harris-rec} Harris ergodic, process
$Y^\ssup{1}$ at a deterministic, sequence of times that increase to
infinity. Still conditioning on $Y^\ssup{2}$, the event that
$Z_{\widetilde{n}_k}\in F$ has probability asymptotically lower
bounded by $\epsilon$. It is then obvious by ergodicity that this
event will occur infinitely often.

According to \cite{A03}, Corollary~VII.3.12, this proves Harris
recurrence of $(Z_n)_{n\in\N}$, and in particular the existence of a
unique invariant measure, \cite{A03}, Theorem~VII.3.5. Now as we want
positive Harris recurrence, we are going to show that this measure is finite.

Note that the previous arguments, together with \cite{A03}, Proposition~VII.3.7, give that
\[
\lim_{N\to\infty}\frac{1}{N}\sum
_{k=1}^N\1_{\{Z_{\widetilde
{n}_k}\in F\}} = {\gamma}
\bigl(F^{\ssup1}\bigr)>0. %
\]

Note that $n_k/k\to2$, since the arrival times of $Y^{\ssup1}$ and
$Y^ {\ssup2}$ are disjoint and have asymptotically the same
distribution. Hence, since $\Wcal\otimes\Wcal\otimes{\Vcal
}(F^\ssup{2})$ is equal to the probability that $Y^{\ssup2}$ hits
$F^{\ssup2}$, we have ${\widetilde{n}_k}/{k}\to2/\Wcal\otimes\Wcal
\otimes{\Vcal}(F^\ssup{2})$ by the ergodic theorem. Noting the
symmetry in the two components, we see that
\[
\lim_{N\to\infty}\frac{1}{N}\sum
_{k=1}^N\1_{\{Z_{k}\in F\}}=\frac{1}2
\bigl( {\gamma}\otimes\Wcal\otimes\Wcal\otimes{\Vcal}(F)+ \Wcal\otimes\Wcal
\otimes{\Vcal}\otimes{\gamma}(F) \bigr). %
\]
Since the right-hand side is a probability measure in $F$, $(Z_n)_{n\in
\N}$ is positive Harris recurrent. Note that we proved the ergodic
theorem in the course of the proof, as well as gave an explicit form
for the invariant measure.

We also see from this proof that the sequence of hitting times of $F$
is non-lattice, since the sequence $(\widetilde n_k)_{k\in\N}$ is
non-lattice, because $( n_k)_{k\in\N}$ is non-lattice.
\end{pf}

\subsection{Longtime average of the connection time}\label{sec-ergodicproofII}

  Here, we give a proof of the ergodic limit in Lemma~\ref
{lem-ergodic} using the Markov chain $Z$ defined in~\eqref{Zdef}. As
we mentioned above, a simpler proof can be done using ergodic theory,
but we will later need the representation of the ergodic limit in terms
of $Z$. We saw in Section~\ref{sec-Mixing} that $(Z_k)_{k\in\N_0}$
is a time-homogeneous, $\psi$-mixing and Harris ergodic Markov chain
on $\Dcal^2$. In this section, we prove the ergodic limit in
Lemma~\ref{lem-ergodic}, giving an explicit formula for the limit
$p_*^{\ssup>}$. The main object in the proof of Theorem~\ref
{thm-longtimebound} in Section~\ref{sec-proofThDeviat} is the
empirical pair measure of $Z$, for which a large-deviation principle is
known to hold.

We\vspace*{1pt} are going to express $\tau_T^{\ssup{\diamond,*}}$ in terms of
$Z$. To this end, we define, for any $z_k= ((x_k^{\ssup1},w_k^{\ssup1},v_k^{\ssup1});(x_k^{\ssup2},w_k^{\ssup2},v_k^{\ssup
2}) )\in\Dcal^2$,
%
\begin{eqnarray}\label{MFdef}
M^{\ssup1}(z_1,z_2)&=&
\frac{\llvert  x_2^{\ssup1}-x_1^{\ssup1}\rrvert  }{v_2^{\ssup
1}},
\nonumber\\[-8pt]\\[-8pt]\nonumber
F_\diamond(z_1,z_2)&=&\int_0^1
\d s \overline\theta^{\ssup
{\diamond}} \bigl(f_*\bigl(p_1(s)\bigr),R
\bigr)\overline\theta^{\ssup{\diamond
}} \bigl(f_*\bigl(p_2(s)\bigr),R
\bigr)\1\bigl\{p_1(s)\mathop{\longleftrightarrow }^{\diamond}_{*}p_2(s)
\bigr\},
\end{eqnarray}
where $p_i(s)=s x_2^{\ssup i}+(1-s)x_1^{\ssup i}$, $s\in[0,1]$,
denotes the path of the $i$th walker from $x_1^{\ssup i}$ to
$x_2^{\ssup i}$. Then $M^{\ssup1}$ is the time that elapses while the
two walkers move from one waypoint arrival to the next one, and
$F_\diamond$ describes the proportion of time that the two are
connected with each other on that way.

Recalling \eqref{Zdef}, we have, for any $n\in\N$,
%
\begin{eqnarray}\label{SNintermsofZ}
S_n=\sum
_{j=1}^n(S_j-S_{j-1})=\sum
_{j=1}^n \frac{\llvert  X^{\ssup
1}_{S_j}-X^{\ssup1}_{S_{j-1}}\rrvert  }{V_{N^{\ssup1}(S_j)}} =\sum
_{j=1}^n M^{\ssup1}(Z_{j-1},Z_j).
\end{eqnarray}
Now we express $\tau_T^{\ssup{\diamond,*}}$ for $T$ replaced by the
waypoint arrival time. For any $n\in\N$, we have
%
\begin{eqnarray}\label{tauintermsofZ}
\tau^{\ssup{\diamond,*}}_{S_n}&=&
\sum_{j=1}^n\int_{S_{j-1}}^{S_j}
\d s \overline\theta^{\ssup{\diamond}}\bigl(f_*\bigl(X^{\ssup1}_s
\bigr),R\bigr)\overline\theta^{\ssup{\diamond
}}\bigl(f_*\bigl(X^{\ssup2}_s
\bigr),R\bigr)\1\bigl\{X^{\ssup1}_s\mathop{\longleftrightarrow
}^{\diamond}_{*}X^{\ssup2}_s\bigr\}\nonumber
\\
&=&\sum_{j=1}^n(S_j-S_{j-1})
\nonumber\\[-8pt]\\[-8pt]\nonumber
&&{}\times
\int_0^1\d s \overline\theta^{\ssup
{\diamond}}
\bigl(f_*\bigl(p_1(s)\bigr),R\bigr)\overline\theta^{\ssup{\diamond}}
\bigl(f_*\bigl(p_2(s)\bigr),R\bigr)\1\bigl\{p_1(s)\mathop{
\longleftrightarrow}^{\diamond
}_{*}p_2(s)\bigr\}
\\
&=&\sum_{j=1}^n M^{\ssup1}(Z_{j-1},Z_j)F_\diamond(Z_{j-1},Z_j),\nonumber
\end{eqnarray}
where $p_i(s)=X^{\ssup i}_{S_{j-1}}+s(X^{\ssup i}_{S_j}-X^{\ssup i}_{S_{j-1}})$.

Now the proof of Lemma~\ref{lem-ergodic} is quite obvious. According
to \cite{A03}, Thoerem~VII.3.6, based on Lemma~\ref{lem-ZHarris},
implies that the distribution of $Z_k$ converges toward its invariant
distribution, which we want to call $\pi$. Hence, $(Z_{j-1},Z_j)$
converges to its invariant distribution $\pi\otimes P$, where we wrote
$P\colon\Dcal\times\Fcal\to[0,1]$ for its transition kernel,
writing $\Fcal$ for the $\sigma$ algebra on $\Dcal$. This
convergence is in total variation sense. Since $M^{\ssup1}$ and
$F_\diamond$ are bounded and measurable, we have that
\[
\lim_{n\to\infty}\frac{1}nS_n=\int
M^{\ssup1} \d(\pi\otimes P)\quad\mbox{and}\quad\lim_{n\to\infty}
\frac{1}n\tau^{\ssup{\diamond,*}}_{S_n}=\int M^{\ssup1}F_\diamond\,
\d(\pi\otimes P). %
\]
Pick $n_T=\sup\{n\in\N\colon S_n\leq T\}$, then it is easy to see
that $\frac{1}Tn_T\to1/\int M^{\ssup1} \,\d(\pi\otimes P)$ as $T\to
\infty$, almost surely and in probability. It is only an exercise to
prove that the above limits are also true if $n$ is replaced by $n_T$.
Furthermore, it is also easy to see that $\frac{1}T(\tau_T^{\ssup
{\diamond,*}}-\tau_{S_{n_T}}^{\ssup{\diamond,*}})$ vanishes almost
surely and in probability as $T\to\infty$. Hence, we have
%
\begin{equation}
\label{ergodlimit} p_*^{\ssup\diamond}=\lim_{T\to\infty}\frac{1}T
\tau_T^{\ssup
{\diamond,*}}=\frac{\int M^{\ssup1}F_\diamond \,\d(\pi\otimes
P)}{\int M^{\ssup1} \,\d(\pi\otimes P)}.
\end{equation}
This ends the proof of Lemma~\ref{lem-ergodic} with the identification
of the limit $p_*^{\ssup\diamond}$ as the right-hand side of \eqref
{ergodlimit}.

\subsection{Proof of Theorem~\texorpdfstring{\protect\ref{thm-longtimebound}}{1.6}}\label
{sec-proofThDeviat}

  Now we turn to the proof of Theorem~\ref{thm-longtimebound},
that is, we prove the upper bound for the downwards deviations of the
normalised connection time, $\frac{1}T\tau_T^{\ssup{>,*}}$, for the
RWP in the limit $T\to\infty$. Let us abbreviate $\tau_T^{\ssup
{>,*}}$ by $\tau_T$. We are going to give an explicit upper bound for
the probability of the event $\{\tau_T\leq T p\}$ for any $p\in
(0,p_*^{\ssup>})$. In order to formulate it, we need to introduce some
more notation, which mostly stems from the theory of large deviations.
See \cite{DemboZeitouni98} for more about this theory.

As a consequence of Lemma~\ref{lem-mixing}, also $(Z_{j-1},Z_j)_{j\in
\N}$ is a $\psi$-mixing and bounded Markov chain. As a nice
consequence, we now have a large-deviation principle (LDP) for the
empirical pair measure of the $Z_n$, defined as
%
\begin{equation}
\label{Qndef} Q_n:=\frac{1}{n} \sum
_{j=1}^n\delta_{(Z_{j-1},Z_j)}\in\mathcal
{M}_1(\Dcal\times\Dcal),
\end{equation}
see \cite{BrycDembo96}, Theorem 1 under the mixing condition (S) and
the remark on page 554, which states that $\psi$-mixing implies (S).
The rate function in \cite{BrycDembo96}, Theorem 1, is given by
\[
I(Q)=\sup_{f\in B(\Dcal^2,\R)} \biggl\{ \int_{\Dcal^2} Q(\d
x,\d y) f(x,y) -\Lambda(f) \biggr\},
\]
where $\Lambda(f)=\lim_{n\to\infty} \frac{1}{n} \log\E_*  [
{\operatorname e }^{\sum_{j=1}^n f(Z_{j-1},Z_j)} ]$, and $B(\Dcal
^2,\R)$ is the set of all bounded, Borel measurable functions on
$\Dcal^2$ to $\R$.

We denote by $\Mcal_1^{\ssup{\rm s}}(\Dcal\times\Dcal)$ the set of
probability measures $Q$ on $\Dcal\times\Dcal$ whose two marginals
coincide. We denote any of the two marginals of such a $Q$ by
$\overline Q$, that is, $\overline Q(A)=Q(A\times\Dcal)=Q(\Dcal
\times A)$ for $A\in\Bcal(\Dcal)$.
Now we use \cite{DemboZeitouni98}, Theorem 6.5.2 for the state space
$\Sigma=\Dcal^2$ and then Theorem 6.5.12 for $k=1$ to identify the
rate function as
%
\begin{equation}
\label{I} I(Q)=H(Q\mid\overline Q\otimes P) =\int_\Dcal
\int_\Dcal Q(\d x,\d y)\log\frac{Q(\d x,\d y)}{
\overline Q(\d x)P(x,\d y)}\qquad\mbox{if }Q
\ll\overline Q\otimes P,
\end{equation}
and $I(Q)=\infty$ otherwise, for $Q\in\Mcal_1^{\ssup{\rm s}}(\Dcal
\times\Dcal)$.

Explicitly, the LDP states that the level sets $\{Q\in\Mcal_1^{\ssup
{\rm s}}(\Dcal\times\Dcal)\colon I(Q)\leq c\}$ are compact for any
$c\in\R$, and that we have the estimates
\[
\limsup_{n\to\infty}\frac{1}n \log\P_*(Q_n\in
F)\leq-\inf_F I\quad\mbox{and}\quad\liminf
_{n\to\infty}\frac{1}n \log\P _*(Q_n\in G)\geq-
\inf_G I, %
\]
for any closed, respectively, open, subset $F$ and $G$ of $\Mcal
_1^{\ssup{\rm s}}(\Dcal\times\Dcal)$.

Theorem~\ref{thm-longtimebound} follows from the following theorem. We
now prefer the notation $\langle f, P\rangle$ for the integral of a
function $f$ with respect to a measure $P$. We recall from \eqref
{ergodlimit} that $p_*^{\ssup>}=\langle M^{\ssup1}F_>,\pi\otimes
P\rangle/\langle M^{\ssup1},\pi\otimes P\rangle$, where $\pi$ is
the invariant distribution of $Z$.

\begin{theorem}\label{thm-longtime} For any $p\in(0,p_*^{\ssup>})$,
%
\begin{equation}
\label{longtime} \limsup_{T\to\infty}\frac{1}T\log\P_*(
\tau_T \leq T p) \leq-\chi_p,
\end{equation}
where
%
\begin{equation}
\label{chidef} \chi_p=\inf\biggl\{\frac{I(Q)}{\langle M^{\ssup1},Q\rangle}\colon Q\in
\Mcal^{\ssup{\rm s}}_1(\Dcal\times\Dcal),\frac{\langle
M^{\ssup1}F_>,Q\rangle}{\langle M^{\ssup1},Q\rangle}\leq p
\biggr\}.
\end{equation}
Moreover, the infimum is attained, and $\chi_p$ is positive.
\end{theorem}

The term $\langle M^{\ssup1},Q\rangle$ is the average time that
elapses between two subsequent arrivals at waypoints, if the two
walkers move in such a way that the distribution of the location,
velocity and next waypoint at two subsequent such arrivals is given by
$Q$, and $\langle M^{\ssup1}F_>,Q\rangle$ is the average portion of
connection time on such a way, and $I(Q)$ is the negative rate of the
probability that the two follow that strategy $Q$ per number of
waypoints. Hence, the upper bound in \eqref{longtime} is intuitive and
can be interpreted. Note that $F_>$ is lower semicontinuous, as the
indicator of connectedness of two points through $\{f_*>\lambda_{\rm
c}(R)\}$ is a countable sum of indicators of open sets. However, in
general $F_>$ may not be upper semicontinuous. This makes it
questionable whether or not also the lower bound in \eqref{longtime}
holds, since the map $Q\mapsto\langle Q,M^{\ssup1} F_>\rangle$ is in
general not continuous.

\begin{pf*}{Proof of Theorem~\ref{thm-longtime}}
That the infimum in \eqref{chidef} is attained is easily seen as
follows. By lower semicontinuity of $F_>$ and a result of Fatou-type
(see, e.g.,~\cite{DemboZeitouni98}, Theorem D.12), the map $Q\mapsto
\langle Q, M^{\ssup1} F_>\rangle$ is also lower semicontinuous. Since
also $I$ is lower semicontinuous and has compact level sets and the map
$Q\mapsto\langle Q, M^{\ssup1}\rangle$ is continuous, it easily
follows that the infimum in \eqref{chidef} is even a minimum.

Now we argue that $\chi_p$ is positive. Indeed, the only minimiser of
$I$ on $\Mcal^{\ssup{\rm s}}_1(\Dcal\times\Dcal)$ is the measure
$\pi\otimes P$, where we recall that $\pi$ is the invariant
distribution of $Z$ and $P$ its transition kernel. To see this, note
that, for any $Q$ satisfying $I(Q)=0$, we have $Q(\d x,\d y)=\overline
Q(\d x)P(x,\d y)$ by the equality discussion in Jensen's inequality,
and from the marginal property it follows that $\overline Q$ is
invariant for $P$, that is, equal to $\pi$ by uniqueness of the
invariant distribution for the chain $Z$. Hence, also the only
minimiser of $Q\mapsto I(Q)/\langle M^{\ssup1},Q\rangle$ is $\pi
\otimes P$, and it satisfies $p_*^{\ssup>}=\langle M^{\ssup1}F_>,\pi
\otimes P\rangle/\langle M^{\ssup1},\pi\otimes P\rangle$, see below
\eqref{ergodlimit}. Therefore, it is not contained in the
admissibility set on the right of \eqref{chidef} and is therefore not
equal to its minimiser. Hence, $\chi_p$ is positive.

Now we prove \eqref{longtime}. We are going to express the time $T$
and the variable $\tau_T$ in terms of integrals over $Q_n$. We write
\[
Z_j= \bigl( \bigl(X^{\ssup1}_{S_j},W_{N^{\ssup1}(S_j)}^{\ssup
1},V_{N^{\ssup1}(S_j)}^{\ssup1}
\bigr), \bigl(X^{\ssup
2}_{S_j},W_{N^{\ssup2}(S_j)}^{\ssup2},V_{N^{\ssup2}(S_j)}^{\ssup
2}
\bigr) \bigr).
\]
%
From \eqref{SNintermsofZ} and \eqref{tauintermsofZ} we have, for any
$n\in\N$,
\[
S_n=n\bigl\langle M^{\ssup1},Q_n\bigr\rangle\quad
\mbox{and}\quad\tau _{S_n}=n\bigl\langle M^{\ssup1}
F_>,Q_n\bigr\rangle, %
\]
recalling the definition of $M^{\ssup1}$ and of $F_>$ in \eqref
{MFdef}, where $p_i(s)=X^{\ssup i}_{S_{j-1}}+s(X^{\ssup
i}_{S_j}-X^{\ssup i}_{S_{j-1}})$.
Hence, we can already give a heuristic proof of Theorem~\ref
{thm-longtime} as follows. The LDP for $(Q_n)_{n\in\N}$ roughly says
that $\P_*(Q_n\approx Q)\approx{\operatorname e }^{-n I(Q)}$ for any
strategy $Q\in\Mcal_1^{\ssup{\rm s}}(\Dcal^2)$. Taking $n$ such
that $T\approx S_n$, we have that $n\approx T/\langle M^{\ssup
1},Q_n\rangle$ and $\tau_T/T\approx\langle M^{\ssup1}
F_>,Q_n\rangle/\langle M^{\ssup1},Q_n\rangle$. Hence, we should have
\begin{eqnarray*}
\P_*(\tau_T\leq pT)&\approx& \P_* \bigl(
\bigl\langle M^{\ssup1} F_>,Q_n\bigr\rangle/\bigl\langle
M^{\ssup1},Q_n\bigr\rangle\leq p \bigr)
\\
&\approx& \exp \biggl(-n\inf\biggl\{I(Q)\colon Q\in\Mcal_1^{\ssup{\rm
s}} \bigl(\Dcal^2\bigr),\frac{\langle M^{\ssup1} F_>,Q\rangle}{\langle
M^{\ssup1},Q\rangle}\leq p\biggr\} \biggr)
\\
&\approx& {\operatorname e }^{-T \chi_p},
\end{eqnarray*}
with $\chi_p$ as defined in Theorem~\ref{thm-longtime}. The main
difficulty in making this line of argument rigorous lies in the
randomness of $n$.

Let us now give a rigorous proof of the upper bound in \eqref
{longtime}. Fix $p\in(0,p_*^{\ssup>})$ and pick a large auxiliary
parameter $K$ and a small one, $\delta>0$. First, we distinguish all
the $n$ no larger than $KT$ such that $T\approx S_n$:
\[
1\leq\sum_{n=\lfloor T/L\rfloor}^{\lfloor KT\rfloor} \1\{S_n
\leq T<S_{n+1}\}+\1\{T\geq S_{\lfloor KT\rfloor+1}\}. %
\]
On the first event, $\{S_n\leq T<S_{n+1}\}$, we have
\[
\tau_T\geq\tau_{S_n}=n\bigl\langle Q_n,M^{\ssup1}F_>
\bigr\rangle\geq (T-L)\frac{\langle M^{\ssup1} F_>,Q_n\rangle}{\langle M^{\ssup
1},Q_n\rangle}\geq T(1-\delta)\frac{\langle M^{\ssup1}
F_>,Q_n\rangle}{\langle M^{\ssup1},Q_n\rangle},
\]
where the last inequality is true for all sufficiently large $T$
(depending only on $\delta$ and $L$), which we want to assume from now.

Observe that $M^{\ssup1}$ is bounded from above by $L=\diam(D)/v_-$,
with probability 1 with respect to $Q$ for any $Q\in\Mcal_1^{\ssup
{\rm s}}(\Dcal^2)$, since $D$ is bounded and all velocities are at
least $v_-$. Hence, we have $S_j-S_{j-1}\leq L$ for any $j\in\N$ and
therefore also $0<S_n/n\leq L$ for any $n\in\N$. Therefore, the
indicator on the event $\{S_n\leq T<S_{n+1}\}$ can be upper bounded as
\[
\1\{S_n\leq T<S_{n+1}\}\leq\1\{T-L\leq S_n\leq
T\}\leq\1\biggl\{ (1-\delta)\frac{T} n\leq\bigl\langle M^{\ssup1},
Q_n\bigr\rangle\leq\frac{T}n\biggr\}. %
\]
This implies the upper bound
\begin{eqnarray*}
\P_*(\tau_T\leq pT)&\leq&\sum_{n=\lfloor T/L\rfloor}^{\lfloor
KT\rfloor}
\P_* \biggl(\frac{\langle M^{\ssup1} F_>,Q_n\rangle
}{\langle M^{\ssup1},Q_n\rangle}\leq\frac{p}{1-\delta},(1-\delta )
\frac{T} n\leq\bigl\langle M^{\ssup1}, Q_n \bigr\rangle
\leq\frac{T}n \biggr)
\\
&&{}+\P_*(T\geq S_{\lfloor KT\rfloor+1}).
\end{eqnarray*}
The last term is an error term, as we will show later that
%
\begin{equation}
\label{errorterm} \lim_{K\to\infty}\limsup_{T\to\infty}
\frac{1}T\log\P_*(T\geq S_{\lfloor KT\rfloor+1})=-\infty.
\end{equation}
Now we cut the sum over $n$ into pieces of length $T\eps$, where $\eps
>0$ is a small auxiliary parameter: 
\[
\sum_{n=\lfloor T/L\rfloor}^{\lfloor KT\rfloor} =\sum
_{i=1+\lfloor
1/L\eps\rfloor}^{\lfloor K/\eps\rfloor} \sum_{(i-1)T\eps<n\leq
iT\eps}.
\]
For fixed $i$ and $(i-1)T\eps<n\leq iT\eps$, we can estimate, for any
large $T$,
%
\begin{eqnarray}\label{LDPappl}
\P_* \biggl(\frac{\langle M^{\ssup1} F_>,Q_n\rangle}{\langle M^{\ssup
1},Q_n\rangle}\leq
\frac{p}{1-\delta},(1-\delta)\frac{T} n\leq \bigl\langle
M^{\ssup1},Q_n\bigr\rangle\leq\frac{T}n \biggr)\leq
\P_*(Q_n\in A_i),
\end{eqnarray}
where
\[
A_i=\biggl\{Q\in\Mcal_1^{\ssup{\rm s}}\bigl(
\Dcal^2\bigr)\colon\frac{\langle
M^{\ssup1} F_>,Q\rangle}{\langle M^{\ssup1},Q\rangle}\leq\frac{p}{1-\delta},
\frac{1-\delta}{i\eps}\leq\bigl\langle M^{\ssup1},Q\bigr\rangle \leq
\frac{1}{(i-1)\eps}\biggr\}. %
\]
Recall that $F_>$ is lower semicontinuous. By \cite{DemboZeitouni98}, Theorem~D.12, the map $Q\mapsto\langle Q,M^{\ssup1}
F_>\rangle$ is also lower semicontinuous. Hence, $A_i$ is closed in
the weak topology. Now we apply the upper bound in the above mentioned
LDP, to obtain, as $T\to\infty$,
\[
\sup_{(i-1)T\eps<n\leq iT\eps}\P_*(Q_n\in A_i)\leq{
\operatorname e }^{-T \widetilde\chi_{p}(\delta,\eps)} {\operatorname e }^{o(T)}, %
\]
where
\begin{eqnarray*}
\widetilde\chi_{p}(\delta,\eps)&=& (i-1)
\eps\inf\bigl\{I(Q)\colon Q\in A_i\bigr\}
\\
&=&(i-1)\eps\inf\biggl\{I(Q)\colon Q\in\Mcal_1^{\ssup{\rm s}}\bigl(
\Dcal ^2\bigr), \frac{\langle M^{\ssup1} F_>,Q\rangle}{\langle M^{\ssup
1},Q\rangle}\leq\frac{p}{1-\delta},
\\
&&{} \frac{1-\delta}{i\eps}\leq\bigl\langle M^{\ssup
1},Q\bigr\rangle\leq\frac{1}{(i-1)\eps}\biggr\}
\\
&\geq& \inf\biggl\{I(Q) \biggl(\frac{1-\delta}{\langle M^{\ssup
1},Q\rangle}-\eps \biggr)\colon Q\in
\Mcal_1^{\ssup{\rm s}}\bigl(\Dcal ^2\bigr),
\frac{\langle M^{\ssup1} F_>,Q\rangle}{\langle M^{\ssup
1},Q\rangle}\leq\frac{p}{1-\delta},
\\
&&{}\frac{1-\delta}{i\eps}\leq\bigl\langle M^{\ssup
1},Q\bigr
\rangle\leq\frac{1}{(i-1)\eps}\biggr\}
\\
&\geq&\inf\biggl\{I(Q) \biggl(\frac{1-\delta}{\langle M^{\ssup
1},Q\rangle}-\eps \biggr)\colon Q\in
\Mcal_1^{\ssup{\rm s}}\bigl(\Dcal ^2\bigr),
\frac{\langle M^{\ssup1} F_>,Q\rangle}{\langle M^{\ssup
1},Q\rangle}\leq\frac{p}{1-\delta}\biggr\}
\\
&=:& \chi_{p}(\delta,\eps).
\end{eqnarray*}
It is easy to see that $\lim_{\eps\downarrow0,\delta\downarrow
0}\chi_{p}(\delta,\eps)=\chi_p$ as defined in \eqref{chidef}.
Hence, the upper bound in~\eqref{longtime} is proved, subject to
\eqref{errorterm}, which we prove now.

Note that\vspace*{1pt} $\{S_n\colon n\in\N_0\}=\{T^{\ssup1}_n\colon n\in\N_0\}
\cup\{T^{\ssup2}_n\colon n\in\N_0\}$, where $T^{\ssup i}_n$ denotes
the arrival time of the $i$th walker at the $n$th waypoint. The $j$th
step $U^{\ssup1}_j$ of the\vspace*{1pt} first of these processes is the duration of
the first walker's travel from the $j$th to the $(j+1)$st waypoint. Hence,
\[
\P_*(T\geq S_{\lfloor KT\rfloor+1})\leq2 \P_*\bigl(T\geq T^{\ssup
1}_{\lfloor KT/2\rfloor+1}
\bigr)\leq2 \P_*\bigl(T\geq\widetilde T^{\ssup
1}_{\lfloor KT/4\rfloor}\bigr),
\]
where\vspace*{-2pt} $\widetilde T^{\ssup1}_n=\sum_{j=0}^{n-1} U^{\ssup1}_{2j}$
denotes the random walk consisting of the even steps only. Hence, we
are now looking at downward deviations of the random walk $(\widetilde
T^{\ssup1}_n)_{n\in\N}$, whose steps $U^{\ssup1}_{2j}$ are
i.i.d.~with support in $[0,L]$. Therefore, Cram\'er's theorem yields
\begin{eqnarray*}
\limsup_{T\to\infty} \frac{1}{T}
\log\P_*\bigl(\widetilde T^{\ssup
1}_{\lfloor KT/4\rfloor}\leq T\bigr) &\leq&
\frac{K}4 \limsup_{T\to\infty} \frac{1}{KT/4} \log\P_*
\biggl(\widetilde T^{\ssup1}_{\lfloor KT/4\rfloor}\leq\frac{4}{K} \lfloor
KT/4\rfloor \biggr)
\\
&\leq&- \frac{K}4 \sup_{\lambda<0} \biggl(\lambda
\frac{4}{K} - \log \E_*\bigl[{\operatorname e }^{\lambda U^{\ssup1}_0}\bigr]
\biggr)
\\
&=&-\sup_{\lambda<0} \biggl(\lambda-\frac{K}4 \log\E_*
\bigl[{\operatorname e }^{\lambda U^{\ssup1}_0}\bigr] \biggr).
\end{eqnarray*}
Note that the essential infimum of $U^{\ssup1}_0$ is equal to zero, as
we assumed that the waypoint measure has a continuous density. Indeed,
if the waypoint walker stands in his waypoint, with probability 1 there
is a non-trivial ball around the location in which the waypoint measure
has a positive density and, therefore, arbitrarily small travels to the
next waypoint have a positive probability.

Hence, $\log\E_*[{\operatorname e }^{\lambda U^{\ssup1}_0}]=o(\llvert
\lambda\rrvert  )$ as $\lambda\to-\infty$ and, therefore, it is possible to
pick a sequence $\lambda_K\to\infty$ as $K\to\infty$ such that
$\lambda_K -\frac{K}4 \log\E_*[{\operatorname e }^{\lambda_K
U^{\ssup1}_0}]\to\infty$ as $K\to\infty$. This implies that \eqref
{errorterm} holds and finishes the proof of Theorem~\ref{thm-longtime}.
\end{pf*}


\section*{Acknowledgements}
Modal'x is affiliated to the LABEX MME-DII.

The authors would like to thank an anonymous referee for his/her
extremely careful reading and the valuable comments to improve the
quality of this article.


%

\printhistory

\begin{thebibliography}{21}
\bibitem{A03}
%
\begin{bbook}[mr]
\bauthor{\bsnm{Asmussen},~\bfnm{S{\o}ren}\binits{S.}}
(\byear{2003}).
\btitle{Applied Probability and Queues},
\bedition{2nd} ed.
\bseries{Applications of Mathematics (New York)}
\bvolume{51}.
\blocation{New York}:
\bpublisher{Springer}.
\bid{mr={1978607}}
\end{bbook}
%

\bptok{imsref}%
\endbibitem

\bibitem{BettstetterWagner02}
%
\begin{barticle}[auto:parserefs-M02]
\bauthor{\bsnm{Bettstetter},~\bfnm{C.}\binits{C.}} \AND
\bauthor{\bsnm{Wagner},~\bfnm{C.}\binits{C.}}
(\byear{2002}).
\btitle{The spatial node distribution of the random waypoint mobility model}.
\bjournal{WMAN},
\bpages{41--58}.
\end{barticle}
%

\bptok{imsref}%
\endbibitem

\bibitem{BettstetterHartensteinCosta04}
%
\begin{barticle}[auto:parserefs-M02]
\bauthor{\bsnm{Bettstetter},~\bfnm{C.}\binits{C.}},
\bauthor{\bsnm{Hartenstein},~\bfnm{H.}\binits{H.}} \AND
\bauthor{\bsnm{P\'erez-Costa},~\bfnm{X.}\binits{X.}} \AND
(\byear{2004}).
\btitle{Stochastic properties of the random waypoint mobility model}.
\bjournal{ACM/Kluwer Wireless Networks 6, Special Issue on Modeling and Analysis of Mobile Networks}
\bvolume{10}
\bpages{555--567}.
\end{barticle}
%

\bptok{imsref}%
\endbibitem

\bibitem{BrycDembo96}
%
\begin{barticle}[mr]
\bauthor{\bsnm{Bryc},~\bfnm{W{\l}odzimierz}\binits{W.}} \AND
\bauthor{\bsnm{Dembo},~\bfnm{Amir}\binits{A.}}
(\byear{1996}).
\btitle{Large deviations and strong mixing}.
\bjournal{Ann. Inst. Henri Poincar\'e Probab. Stat.}
\bvolume{32}
\bpages{549--569}.
\bid{issn={0246-0203}, mr={1411271}}
\end{barticle}
%

\bptok{imsref}%
\endbibitem

\bibitem{CBD02}
%
\begin{barticle}[auto:parserefs-M02]
\bauthor{\bsnm{Camp},~\bfnm{T.}\binits{T.}},
\bauthor{\bsnm{Boleng},~\bfnm{J.}\binits{J.}} \AND
\bauthor{\bsnm{Davies},~\bfnm{V.}\binits{V.}}
(\byear{2002}).
\btitle{A survey of mobility models for ad-hoc network research}
\bjournal{WCMC: Special Issue on Mobile ad-hoc Networking: Research,
Trends and Applications}
\bvolume{2}
\bpages{483--502}.
\end{barticle}
%

\bptok{imsref}%
\endbibitem

\bibitem{Clementi2009}
%
\begin{bincollection}[mr]
\bauthor{\bsnm{Clementi},~\bfnm{Andrea~E.~F.}\binits{A.E.F.}},
\bauthor{\bsnm{Pasquale},~\bfnm{Francesco}\binits{F.}} \AND
\bauthor{\bsnm{Silvestri},~\bfnm{Riccardo}\binits{R.}}
(\byear{2009}).
\btitle{M{ANETS}: High mobility can make up for low transmission power}.
In \bbooktitle{Automata, Languages and Programming. {P}art {II}}.
\bseries{Lecture Notes in Computer Science}
\bvolume{5556}
\bpages{387--398}.
\blocation{Berlin}:
\bpublisher{Springer}.
\bid{doi={10.1007/978-3-642-02930-1_32}, mr={2544811}}
\end{bincollection}
%

\bptok{imsref}%
\endbibitem

\bibitem{DemboZeitouni98}
%
\begin{bbook}[mr]
\bauthor{\bsnm{Dembo},~\bfnm{Amir}\binits{A.}} \AND
\bauthor{\bsnm{Zeitouni},~\bfnm{Ofer}\binits{O.}}
(\byear{2010}).
\btitle{Large Deviations Techniques and Applications}.
\bseries{Stochastic Modelling and Applied Probability}
\bvolume{38}.
\blocation{Berlin}:
\bpublisher{Springer}.
\bnote{Corrected reprint of the second (1998) edition}.
\bid{doi={10.1007/978-3-642-03311-7}, mr={2571413}}
\end{bbook}
%

\bptok{imsref}%
\endbibitem

\bibitem{HyytiaLassilaVirtamo06}
%
\begin{barticle}[auto:parserefs-M02]
\bauthor{\bsnm{Hyytia},~\bfnm{E.}\binits{E.}},
\bauthor{\bsnm{Lassila},~\bfnm{P.}\binits{P.}} \AND
\bauthor{\bsnm{Virtamo},~\bfnm{J.}\binits{J.}}
(\byear{2006}).
\btitle{Spatial node distribution of the random waypoint mobility
model with applications}.
\bjournal{IEEE Trans. Mob. Comput.}
\bvolume{5}
\bpages{680--694}.
\end{barticle}
%

\bptok{imsref}%
\endbibitem

\bibitem{KaspiMandelbaum94}
%
\begin{barticle}[mr]
\bauthor{\bsnm{Kaspi},~\bfnm{Haya}\binits{H.}} \AND
\bauthor{\bsnm{Mandelbaum},~\bfnm{Avi}\binits{A.}}
(\byear{1994}).
\btitle{On {H}arris recurrence in continuous time}.
\bjournal{Math. Oper. Res.}
\bvolume{19}
\bpages{211--222}.
\bid{doi={10.1287/moor.19.1.211}, issn={0364-765X}, mr={1290020}}
\end{barticle}
%

\bptok{imsref}%
\endbibitem

\bibitem{Kesten74}
%
\begin{barticle}[mr]
\bauthor{\bsnm{Kesten},~\bfnm{Harry}\binits{H.}}
(\byear{1974}).
\btitle{Renewal theory for functionals of a {M}arkov chain with
general state space}.
\bjournal{Ann. Probab.}
\bvolume{2}
\bpages{355--386}.
\bid{mr={0365740}}
\end{barticle}
%

\bptok{imsref}%
\endbibitem

%

\bibitem{Boudec04}
%
\begin{barticle}[auto:parserefs-M02]
\bauthor{\bsnm{Le Boudec},~\bfnm{J.-Y.}\binits{J.-Y.}}
(\byear{2007}).
\btitle{Understanding the simulation of mobility models with palm calculus}.
\bjournal{Perform. Eval.}
\bvolume{64}
\bpages{126--146}.
\end{barticle}
%

\bptok{imsref}%
\endbibitem

\bibitem{BoudecVojnovic06}
%
\begin{barticle}[auto:parserefs-M02]
\bauthor{\bsnm{Le Boudec},~\bfnm{J.-Y.}\binits{J.-Y.}} \AND
\bauthor{\bsnm{Vojnovic},~\bfnm{M.}\binits{M.}}
(\byear{2006}).
\btitle{The random trip model: Stability, stationary regime, and
perfect simulation}.
\bjournal{IEEE/ACM Transactions on Networking}
\bvolume{14}
\bpages{1153--1166}.
\end{barticle}
%

\bptok{imsref}%
\endbibitem

\bibitem{MeesterRoy96}
%
\begin{bbook}[mr]
\bauthor{\bsnm{Meester},~\bfnm{Ronald}\binits{R.}} \AND
\bauthor{\bsnm{Roy},~\bfnm{Rahul}\binits{R.}}
(\byear{1996}).
\btitle{Continuum Percolation}.
\bseries{Cambridge Tracts in Mathematics}
\bvolume{119}.
\blocation{Cambridge}:
\bpublisher{Cambridge Univ. Press}.
\bid{doi={10.1017/CBO9780511895357}, mr={1409145}}
\end{bbook}
%

\bptok{imsref}%
\endbibitem

\bibitem{Penrose03}
%
\begin{bbook}[mr]
\bauthor{\bsnm{Penrose},~\bfnm{Mathew D.}\binits{M.D.}}
(\byear{2003}).
\btitle{Random Geometric Graphs}.
\bseries{Oxford Studies in Probability}
\bvolume{5}.
\blocation{Oxford}:
\bpublisher{Oxford Univ. Press}.
\bid{doi={10.1093/acprof:oso/9780198506263.001.0001}, mr={1986198}}
\end{bbook}
%

\bptok{imsref}%
\endbibitem

\bibitem{Penrose91}
%
\begin{barticle}[mr]
\bauthor{\bsnm{Penrose},~\bfnm{Mathew~D.}\binits{M.D.}}
(\byear{1991}).
\btitle{On a continuum percolation model}.
\bjournal{Adv. in Appl. Probab.}
\bvolume{23}
\bpages{536--556}.
\bid{doi={10.2307/1427621}, issn={0001-8678}, mr={1122874}}
\end{barticle}
%

\bptok{imsref}%
\endbibitem

\bibitem{Penrose95}
%
\begin{barticle}[mr]
\bauthor{\bsnm{Penrose},~\bfnm{Mathew~D.}\binits{M.D.}}
(\byear{1995}).
\btitle{Single linkage clustering and continuum percolation}.
\bjournal{J. Multivariate Anal.}
\bvolume{53}
\bpages{94--109}.
\bid{doi={10.1006/jmva.1995.1026}, issn={0047-259X}, mr={1333129}}
\end{barticle}
%

\bptok{imsref}%
\endbibitem

\bibitem{Peres13}
%
\begin{barticle}[mr]
\bauthor{\bsnm{Peres},~\bfnm{Yuval}\binits{Y.}},
\bauthor{\bsnm{Sinclair},~\bfnm{Alistair}\binits{A.}},
\bauthor{\bsnm{Sousi},~\bfnm{Perla}\binits{P.}} \AND
\bauthor{\bsnm{Stauffer},~\bfnm{Alexandre}\binits{A.}}
(\byear{2013}).
\btitle{Mobile geometric graphs: Detection, coverage and percolation}.
\bjournal{Probab. Theory Related Fields}
\bvolume{156}
\bpages{273--305}.
\bid{doi={10.1007/s00440-012-0428-1}, issn={0178-8051}, mr={3055260}}
\end{barticle}
%

\bptok{imsref}%
\endbibitem

\bibitem{QZ07}
%
\begin{barticle}[auto:parserefs-M02]
\bauthor{\bsnm{Quintanilla},~\bfnm{J.~A.}\binits{J.A.}} \AND
\bauthor{\bsnm{Ziff},~\bfnm{R.~M.}\binits{R.M.}}
(\byear{2007}).
\btitle{Asymmetry in the percolation thresholds of fully penetrable
disks with two different radii}.
\bjournal{Phys. Rev. E}
\bvolume{76}
\bpages{051115}.
\end{barticle}
%

\bptok{imsref}%
\endbibitem

\bibitem{Roy11}
%
\begin{bbook}[auto:parserefs-M02]
\bauthor{\bsnm{Roy},~\bfnm{R.~R.}\binits{R.R.}}
(\byear{2011}).
\btitle{Handbook of Mobile Ad Hoc Networks for Mobility Models}.
\blocation{New York}:
\bpublisher{Springer}.
\end{bbook}
%

\bptok{imsref}%
\endbibitem


\bibitem{Sarkar97}
%
\begin{barticle}[mr]
\bauthor{\bsnm{Sarkar},~\bfnm{Anish}\binits{A.}}
(\byear{1997}).
\btitle{Continuity and convergence of the percolation function in
continuum percolation}.
\bjournal{J. Appl. Probab.}
\bvolume{34}
\bpages{363--371}.
\bid{issn={0021-9002}, mr={1447341}}
\end{barticle}
%

\bptok{imsref}%
\endbibitem
\end{thebibliography}
\end{document}